\documentclass[aop,seceqn,citesort,MSNbibl,dvips]{arximspdf}
\usepackage{graphics}
\doi{10.1214/10-AOP543}
\volume{39}
\issue{1}
\pubyear{2011}
\firstpage{252}
\lastpage{290}

\makeatletter

\newproclaim{convention}{Convention}[section]
\newproclaim{remark}[convention]{Remark}
\newproclaim{example}{Example}
\newtheorem{proposition}[convention]{Proposition}
\newproclaim{definition}[convention]{Definition}
\newtheorem{theorem}[convention]{Theorem}
\newtheorem{conj}[convention]{Conjecture}
\newtheorem{lemma}[convention]{Lemma}

\newcommand{\R}{\mathbb{R}}
\newcommand{\N}{\mathbb{N}}
\newcommand{\Z}{\mathbb{Z}}
\newcommand{\E}{\mathbb{E}}
\renewcommand{\P}{\mathbb{P}}
\newcommand{\M}{{\mathcal M}}
\newcommand{\one}{\mathbf{1}}
\newcommand{\bu}{\mathbf{u}}
\newcommand{\bv}{\mathbf{v}}

\makeatother

\begin{document}
\begin{frontmatter}

\title{On the moments and the interface of the symbiotic branching model}
\runtitle{On the symbiotic branching model}

\begin{aug}
\author[A]{\fnms{Jochen} \snm{Blath}\ead[label=e1]{blath@math.tu-berlin.de}},
\author[A]{\fnms{Leif} \snm{D\"oring}\corref{}\thanksref{t2}\ead[label=e2]{doering@math.tu-berlin.de}} and
\author[B]{\fnms{Alison} \snm{Etheridge}\ead[label=e3]{etheridg@stats.ox.ac.uk}}
\runauthor{J. Blath, L. D\"oring and A. Etheridge}
\affiliation{Technische Universit\"at Berlin, Technische Universit\"at
Berlin\break
and University of Oxford}
\address[A]{J. Blath\\
L. D\"oring\\
Institut f\"ur Mathematik\\
Technische Universit\"at Berlin\\
Stra\ss e des 17 Juni 136\\
10623 Berlin\\
Germany\\
\printead{e1}\\
\phantom{E-mail: }\printead*{e2}}
\address[B]{A. Etheridge\\
Department of Statistics\\
University of Oxford\\
1, South Parks Road\\
Oxford OX1 3TG\\
United Kingdom\\
\printead{e3}}
\end{aug}

\thankstext{t2}{Supported in part by the DFG International Research
Training Group
``Stochastic Models of Complex Processes'' and the Berlin Mathematical School.}

\received{\smonth{7} \syear{2009}}
\revised{\smonth{2} \syear{2010}}

%
\begin{abstract}
In this paper we introduce a critical curve separating the asymptotic
behavior of the moments of the symbiotic branching model, introduced by
Etheridge and Fleischmann [\textit{Stochastic Process. Appl.}
\textbf{114} (2004) 127--160] into two regimes. Using arguments based
on two
different dualities and a classical result of Spitzer [\textit{Trans. Amer.
Math. Soc.} \textbf{87} (1958) 187--197] on the exit-time of a
planar Brownian motion from a wedge, we prove that the parameter
governing the model provides regimes of bounded and exponentially
growing moments separated by subexponential growth. The moments turn
out to be closely linked to the limiting distribution as time tends to
infinity. The limiting distribution can be derived by a self-duality
argument extending a result of Dawson and Perkins [\textit{Ann. Probab.}
\textbf{26} (1998) 1088--1138] for the mutually catalytic
branching model.

As an application, we show how a bound on the $35$th moment improves
the result of Etheridge and Fleischmann [\textit{Stochastic
Process. Appl.} \textbf{114} (2004) 127--160] on the speed of the propagation
of the interface of the symbiotic branching model.
\end{abstract}

%
\begin{keyword}[class=AMS]
\kwd[Primary ]{60K35}
\kwd[; secondary ]{60J80}.
\end{keyword}
\begin{keyword}
\kwd{Symbiotic branching model}
\kwd{mutually catalytic branching}
\kwd{stepping stone model}
\kwd{parabolic Anderson model}
\kwd{moment duality}
\kwd{self-duality}
\kwd{propagation of interface}
\kwd{exit distribution}.
\end{keyword}

\end{frontmatter}

\section{Introduction}\label{sec:intro}

In 2004, Etheridge and Fleischmann \cite{EF04} introduced a stochastic
spatial model of two interacting populations known as the symbiotic
branching model, parametrized by a parameter $\varrho\in[-1,1]$
governing the correlation between the two driving noises. The model can
be considered in three different spatial setups which we now explain.

First, the continuous-space symbiotic branching model is given by the
system of stochastic partial differential equations
%
%
\begin{equation}\label{eqn:spde}\quad
{\mathrm{cSBM}(\varrho,\kappa)}_{u_0, v_0}\dvtx
\cases{
\dfrac{\partial}{\partial t}u_t(x) = \dfrac{1}{2}\Delta u_t(x) +
\sqrt{ \kappa u_t(x) v_t(x)} \,dW^1_t(x),\vspace*{2pt}\cr
\dfrac{\partial}{\partial t}v_t(x) = \dfrac{1}{2}\Delta v_t(x) +
\sqrt{ \kappa u_t(x) v_t(x)} \,dW^2_t(x),\vspace*{2pt}\cr
u_0(x) \ge0, \qquad x \in\R,\cr
v_0(x) \ge0, \qquad\hspace*{1.1pt} x \in\R,}
\end{equation}
where $\Delta$ denotes the Laplace operator and $\kappa> 0$ is a fixed
constant known as the branching
rate. $\mathbf{{W}} = ({W}^1,
{W}^2)$ is a pair of correlated standard Gaussian white noises on
$\mathbb{R}_+ \times\mathbb{R}$ with correlation $\varrho
\in[-1,1]$, that is, the unique Gaussian process with covariance structure
%
%
\begin{eqnarray}\label{correlation}
\mathbb{E} [ W^1_{t_1}(A_1)W^1_{t_2}(A_2) ] &=& (t_1\wedge t_2)
\ell(A_1\cap A_2),\\
\mathbb{E} [ W^2_{t_1}(A_1)W^2_{t_2}(A_2) ] &=& (t_1\wedge t_2)
\ell(A_1\cap A_2),\\
\mathbb{E} [ W^1_{t_1}(A_1)W^2_{t_2}(A_2) ] &=& \varrho
(t_1\wedge t_2) \ell(A_1\cap A_2),
\end{eqnarray}
where $\ell$ denotes Lebesgue measure, $A_1,A_2\in\mathcal{B}(\R)$
and $t_1,t_2 \geq0$.
Note that we work with a white noise $\mathbf{{W}}$ in the sense of
Walsh \cite{W86}. Solutions of this model have been considered
rigorously in the framework of the corresponding martingale problem in
Theorem 4 of \cite{EF04}, which states that, under suitable conditions
on the initial conditions $u_0(\cdot), v_0(\cdot)$, a solution exists
for all $\varrho\in[-1, 1]$. The martingale problem is well posed for
all $\varrho\in[-1,1)$, which implies the strong Markov property
except in the boundary case $\varrho=1$.

For a discrete spatial version we consider the system of interacting
diffusions on $\Z^d$, with values in $\R_{\geq0}$, defined by the
coupled stochastic differential equations
%
%
\begin{equation}\label{eq:dsbm}
\mathrm{dSBM}(\varrho,\kappa)_{u_0, v_0} \dvtx
\cases{
d u_t(i)=\Delta u_t(i) \,dt + \sqrt{\kappa u_t(i)v_t(i)}
\,dB^1_t(i),\vspace*{2pt}\cr
d v_t(i)=\Delta v_t(i) \,dt + \sqrt{\kappa
u_t(i)v_t(i)} \,dB^2_t(i),\vspace*{2pt}\cr
u_0(i) \ge0, \qquad i \in\Z^d,\cr
v_0(i) \ge0, \qquad\hspace*{1.1pt} i \in\Z^d,}
\end{equation}
where now $ \{B^1(i),B^2(i) \}_{i \in\Z^d}$ is a family of
standard Brownian motions with covariances given by
%
%
\begin{equation}\label{eq:dcv}
[ B_{\cdot}^n(i),B_{\cdot}^m(j)]_t =
\cases{
\varrho t, &\quad $i=j$ and $n \neq m$,\cr
t, &\quad $i=j$ and $n=m$, \cr
0, &\quad otherwise.}
\end{equation}
In the discrete case, $\Delta$ denotes the discrete Laplacian
\[
\Delta u_t(i)=\sum_{|k-i|=1}\frac{1}{2d}\bigl(u_t(k)-u_t(i)\bigr).
\]

Note that in this paper we denote by $[N_{\cdot},M_{\cdot}]_t$ the
cross-variation of two martingales $N,M$. This is to avoid confusion
with $\langle f,g\rangle$ which will be defined to be the sum (resp.,
integral) of the product of $f$ and $g$.

Finally, the nonspatial symbiotic branching model is defined by the
stochastic differential equations
\[
\mathrm{SBM}(\varrho,\kappa)_{u_0, v_0} \dvtx
\cases{
du_t=\sqrt{\kappa u_tv_t} \,dB^1_t,\cr
dv_t=\sqrt{\kappa u_tv_t} \,dB^2_t,\cr
u_0\geq0,\cr
v_0\geq0.}
\]
Again, the noises are correlated with $[B^1_{\cdot},B^2_{\cdot
}]_t=\varrho t$. This simple toy-model (see also \cite{R95} and \cite
{DFX05}) can be analyzed quite simply and will be used to prove
properties of the spatial models.
\begin{convention}
From time to time we skip the dependence on $\varrho, \kappa, u_0$ and
$v_0$ if there is no ambiguity. Solutions of $\mathrm{cSBM}, \mathrm
{SBM}$ and $\mathrm{dSBM}$ for $d\leq2$ are called symbiotic branching
processes in the recurrent case whereas solutions of $\mathrm{dSBM}$
for $d\geq3$ are called symbiotic branching processes in the transient case.
\end{convention}

Interestingly, symbiotic branching models include well-known spatial
models from different branches of probability theory. In the discrete
spatial case (and analogously in continuous-space) interacting
diffusions of the type
%
%
\begin{equation}
\label{int}
d w_t(i)=\Delta w_t(i) \,dt + \sqrt{\kappa f(w_t(i))} \,dB_t(i)
\end{equation}
have been studied extensively in the literature. Some important
examples are the following:
\begin{example}\label{ex1}
The stepping stone model from mathematical genetics: $f(x)= x(1-x)$.
\end{example}
\begin{example}\label{ex2}
The parabolic Anderson model (with Brownian potential) from
mathematical physics: $f(x)= x^2$.
\end{example}
\begin{example} \label{ex3}
The super random walk from probability theory: $f(x)= x$.
\end{example}

For the super random walk, $\kappa$ is the branching rate which in this
case is time--space independent. In \cite{DP98}, a two-type model based
on two super random walks with time--space dependent branching was
introduced. The branching rate for one species is proportional to the
value of the other species. More precisely, the authors considered
\begin{eqnarray*}
d u_t(i)&=&\Delta u_t(i) \,dt + \sqrt{\kappa u_t(i)v_t(i)}\,
dB^1_t(i),\\
d v_t(i)&=&\Delta v_t(i) \,dt + \sqrt{\kappa u_t(i)v_t(i)}\, dB^2_t(i),
\end{eqnarray*}
where now $ \{B^1(i),B^2(i) \}_{i \in\Z^d}$ is a family of
independent standard Brownian motions. Solutions are called mutually
catalytic branching processes. In the following years, properties of
this model were well studied (see, e.g., \cite{CK00} and \cite
{CDG04}). The corresponding continuous-space version was also treated
in \cite{DP98}.

For correlation $\varrho=0$, solutions of the symbiotic branching model
are obviously solutions of the mutually catalytic branching model. The
case $\varrho=-1$ with the additional assumption $u_0+v_0\equiv1$
corresponds to the stepping stone model. To see this, observe that in
the perfectly negatively correlated case $B^1(i)=-B^2(i)$ which implies
that the sum $u+v$ solves a discrete heat equation and with the further
assumption $u_0+v_0\equiv1$ stays constant for all time. Hence, for
all $t\geq0$, $u(t,\cdot)\equiv1-v(t,\cdot)$, which shows that $u$ is
a solution of the stepping stone model with initial condition $u_0$ and
$v$ is a solution with initial condition $v_0$. Finally, suppose $w$ is
a solution of the parabolic Anderson model, then, for $\varrho=1$, the
pair $(u,v):=(w,w)$ is a solution of the symbiotic branching model with
initial conditions \mbox{$u_0=v_0=w_0$}.

The purpose of this and the accompanying paper \cite{AD09} is to
understand the nature of the symbiotic branching model better. How does
the model depend on the correlation $\varrho$? Are properties of the
extremal cases $\varrho\in\{-1,0,1\}$ inherited by some parts of the
parameter space? Since the longtime behavior of the super random walk,
stepping stone model, mutually catalytic branching model and parabolic
Anderson model is very different, one might guess that the parameter
space $[-1,1]$ can be divided into disjoint subsets corresponding to
different regimes.

The focus of \cite{AD09} is second moment properties. In the discrete
setting, but with a more general setup, growth of second moments is
analyzed in detail. A moment duality is used to reduce the problem to
moment generating functions and Laplace transforms of local times of
discrete-space Markov processes. A precise analysis of those is used to
derive intermittency and aging results which show that different
regimes occur for $\varrho<0$, $\varrho=0$ and $\varrho>0$.

In contrast to \cite{AD09}, the present paper is not restricted to
second moment properties. The aim is to understand the pathwise
behavior of symbiotic branching processes better.
\begin{remark}\label{rem:migration}
In this paper, we restrict ourselves to the simplest setups which
already provide the full variety of results. For the discrete spatial
model we thus restrict ourselves to the discrete Laplacian instead of
allowing more general transitions. This is not necessary; see \cite
{DP98} or \cite{CDG04} for a construction of solutions and main
properties for more general underlying migration mechanisms in the case
$\varrho=0$. Furthermore, we mainly restrict ourselves to homogeneous
initial conditions and remark where results hold more generally. Here,
for nonnegative real numbers we denote by $\bu$ the constant functions
$u(\cdot)\equiv u$.
\end{remark}

The paper is organized as follows: our main results are presented in
Section \ref{subsec:smr}. Before proving the results, we collect basic
properties of the symbiotic branching models and discuss the dualities
that we need. This is carried out in Section \ref{sec:bpd}. The final
sections are devoted to the proofs. In Section \ref{sec:comnvlaw},
proofs of the longtime convergence in law are given, and in Section
\ref{sec:moments} we discuss the longtime behavior of moments. Finally, in
Section \ref{sec:wavespeed} we show how to use the results of
Section \ref{sec:moments} to strengthen the main result of \cite{EF04}.

\section{Results}\label{subsec:smr}
Before stating the main results, we briefly recall from \cite{EF04}
that the state space of $\mathrm{cSBM}$ is given by pairs of tempered
functions, that is, pairs of functions contained in
\[
M_{\mathrm{tem}}= \Bigl\{u | u\dvtx\R\to\R_{\geq0}, \lim_{|x|\rightarrow\infty
}u(x)\phi_{\lambda}(x)\mbox{ exists and }\Vert u\phi_{\lambda
}\Vert
_{\infty
}<\infty\ \forall\lambda<0 \Bigr\},
\]
where $\phi_{\lambda}(x)=e^{\lambda|x|}$, and we think of $M_{\mathrm{tem}}$ as
being topologized by the metric given in \cite{EF04}, equation (13),
yielding a Polish space.

The state space for $\mathrm{dSBM}$ is similar. It was not discussed in
\cite{EF04} and so we present details in Section \ref{sec:bpd}.


\subsection{Convergence in law}

We begin with a result, generalizing Theorem 1.5 of~\cite{DP98}, on the
longtime behavior of the laws of symbiotic branching processes in the
recurrent case.
\begin{proposition}\label{prop:convlaw}
Suppose $(u_t,v_t)$ is a spatial symbiotic branching process in the
recurrent case with $\varrho\in(-1,1)$, $\kappa>0$ and initial
conditions $u_0=\bu, v_0=\bv$. Let $B^1$ and $B^2$ be two Brownian
motions with covariance
\[
[B^1_{\cdot},B^2_{\cdot}]_t=\varrho t,\qquad t\ge0,
\]
and initial conditions $B_0^1=u, B_0^2=v$. Further, let
\[
\tau=\inf\{t\geq0\dvtx B^1_t B^2_t=0 \}
\]
be the first exit time of the correlated Brownian motions $B^1, B^2$
from the upper right quadrant.
Then, weakly in $M_{\mathrm{tem}}^2$,
\[
\P^{\bu,\bv}[(u_t,v_t)\in\cdot]\Rightarrow P^{u,v}[(\bar B^1_{\tau
},\bar B^2_{\tau})\in\cdot]
\]
as $t\to\infty$. Here, $(\bar B^1_\tau, \bar B^2_\tau)$ denotes the
pair of constant functions on $\R$, respectively, $\Z^d$ ($d=1,2$)
taking the
values of the stopped Brownian motions $(B^1_\tau, B^2_\tau)$.
\end{proposition}

In particular, the proposition shows ultimate extinction of one species
in law.
\begin{remark}\label{all}
For simplicity, Proposition \ref{prop:convlaw} is formulated for
constant initial conditions even though the result holds more
generally. Theorem 1.5 of \cite{DP98} (the case $\varrho=0$) was
extended in \cite{CKP00} to nondeterministic initial conditions: for
fixed $u,v\geq0$ let $\M_{u,v}$ be the set of probability measures
$\nu
$ on $M_{\mathrm{tem}}^2$ such that
%
%
\begin{equation}
\sup_{x \in\R} \int\bigl(a^2(x) + b^2(x)\bigr) \,d\nu(a,b)<\infty
\end{equation}
and
%
%
\begin{equation}\qquad
\lim_{t \to\infty} \int\bigl[ \bigl(P_ta(x)- u\bigr)^2 + \bigl(P_tb(x)-v\bigr)^2\bigr]\, d\nu
(a,b) =0 \qquad\mbox{for all } x \in\R.
\end{equation}
Here, $(P_t)$ denotes the transition semigroup of Brownian motion (the
definition for the discrete case is similar). The proof of \cite{CKP00}
can also be applied to $\varrho\neq0$ and, thus, Proposition \ref
{prop:convlaw} holds in the same way for initial distributions $\nu\in
\mathcal{M}_{u,v}$.
\end{remark}

The restriction to $\varrho\in(-1,1)$ arises from our method of proof
which exploits a self-duality of the process which gives no information
for $\varrho\in\{-1,1\}$. Let us briefly discuss the behavior of the
limiting distributions in the boundary cases $\varrho\in\{-1, 1\}$
which are well known in the literature and fit neatly into our result.
First, suppose $(w_t)$ is a solution of the stepping stone model (see
Example \ref{ex1}) and $w_0\equiv w \in[0,1]$. It was proved in
\cite{S80} that
%
%
\begin{equation}
\label{eq:shigalimit}
\mathcal{L}^{\mathbf{w}}(w_t)\stackrel{t\rightarrow\infty
}{\Rightarrow} w
\delta_{\one}
+(1-w) \delta_{\mathbf{0}},
\end{equation}
where $\delta_{\one}$ (resp., $\delta_{\mathbf{0}}$) denotes the Dirac
distribution concentrated on the constant function $\one$ (resp.,
$\mathbf0$). This can be reformulated in terms of perfectly anti-correlated
Brownian motions $(B^1, B^2)$ as before: for $\varrho=-1$, the pair
$(B^1, B^2)$ takes values only on the straight line connecting $(0,1)$
and $(1,0)$, and stops at the boundaries. Hence, the law of $(B^1_{\tau
}, B^2_{\tau})$ is a mixture of $\delta_{(0,1)}$ and $\delta_{(1,0)}$
and the probability of hitting $(1,0)$ is equal to the probability of a
one-dimensional Brownian motion started in $w \in[0,1]$ hitting $1$
before $0$, which is $w$, and hence matches (\ref{eq:shigalimit}).
Second, let $(w_t)$ be a solution of the parabolic Anderson model with
Brownian potential (see Example \ref{ex2}) and constant initial
condition $w_0 \equiv w \geq0$. In \cite{r7} it was shown that
%
%
\begin{equation}
\mathcal{L}^{\mathbf{w}}(w_t)\stackrel{t\rightarrow\infty
}{\Rightarrow}
\delta_{\mathbf{0}}.
\end{equation}
As discussed above, when viewed as a symbiotic branching process with
$\varrho=1$, this implies
%
%
\begin{equation}
\mathcal{L}^{\mathbf w, {\mathbf w}}(u_t,v_t)\stackrel{t\rightarrow
\infty
}{\Rightarrow} \delta_{\mathbf{0}, \mathbf{0}}.
\end{equation}
From the viewpoint of two perfectly positive-correlated Brownian
motions, we obtain the same result since they simple move on the
diagonal dissecting the upper right quadrant until they eventually get
absorbed in the origin, that is, $(B^1_\tau, B^2_\tau)=(0,0)$ almost surely.

To summarize, we have seen that the weak longtime behavior (in the
recurrent case) of the classical models connected to symbiotic
branching is appropriately described by correlated Brownian motions
hitting the boundary of the upper right quadrant.

\subsection{Nonalmost-sure behavior}
In contrast to extinction in law, the almost-sure behavior is very
different. In the recurrent case for the mutually catalytic branching
model, Cox and Klenke \cite{CK00} showed that, almost surely, there is
no longtime local extinction of any type, but in fact the locally
predominant type changes infinitely often. It is not hard to see that
the same is true for symbiotic branching with $\varrho\in(-1,1)$. We do
not give a proof since it follows from Proposition \ref{prop:convlaw}
along the same lines as in \cite{CK00}.
\begin{proposition}\label{prop:pb}
Let $\varrho\in(-1,1)$, $\kappa>0$ and suppose $(u_t,v_t)$ is a
spatial symbiotic branching process in the recurrent case with initial
distribution $u_0=\bu, v_0=\bv$. Then, for all $(u', v') \in\{(x, 0)
\dvtx x\in\R_{\geq0}\} \cup\{ (0,y) \dvtx y\in\R_{\geq0}\}$ and
$K \subset\R$ bounded,
\[
\P^{\bu,\bv} \Bigl[{\liminf_{t\rightarrow\infty} \sup_{x\in K}}
\Vert(u_t(x), v_t(x))-(u', v')\Vert=0 \Bigr] =1,
\]
respectively, for $K\subset\Z^d$ bounded,
\[
\P^{\bu,\bv} \Bigl[ {\liminf_{t\rightarrow\infty} \sup_{k\in K}}
\Vert(u_t(k), v_t(k))-(u', v')\Vert=0 \Bigr] =1.
\]
\end{proposition}

Again, as in Remark \ref{all}, the result holds for random initial
conditions of the class $\mathcal{M}_{u,v}$. Note that Proposition
\ref
{prop:pb} depends strongly on the spatial structure since in the
nonspatial model almost sure convergence holds (see Proposition \ref
{prop:sconv}).

\subsection{Longtime behavior of moments}
In \cite{AD09} the second moments of symbiotic branching processes are
analyzed. This particular case admits a detailed study since a moment
duality (see Lemma \ref{la:mdual}) has a particularly simple structure
which allows one to reduce the study of the moments to that of moment
generating functions and Laplace transforms of local times. Here we are
interested in the behavior of moments as $t$ tends to infinity. The two
available dualities (self-duality and moment duality) are combined in
two steps. First, a self-duality argument combined with an equivalence
between bounded moments of the exit time distribution and of the exit
point distribution for correlated Brownian motions stopped on exiting
the first quadrant is used to understand the effect of $\varrho$. It
turns out that for any $p>1$ there are critical values, independent of
$\kappa$, dividing regimes in which the moments $\E^{1,1}[u_t^p]$,
$\E
^{\one,\one}[u_t(k)^p]$ and $\E^{\one,\one}[u_t(x)^p]$ are bounded
in $t$ or
grow to infinity. Second, for $p\in\N$, a perturbation argument
combined with the first step and a moment duality is used to analyze
the growth to infinity in more detail.

The following critical curve captures the effect of $\varrho$. Note
that the definition is independent of $\kappa$ which will become
important in the second step.
\begin{definition} \label{def:cc}
We define the critical curve of symbiotic branching models to be the
real-valued function $p\dvtx(-1,1) \to\R^+$, given by
%
%
\begin{equation}\label{criticalcurve}
p(\varrho)=\frac{\pi}{{\pi}/{2}+\arctan({\varrho
}/({\sqrt
{1-\varrho^2}}))}.
\end{equation}
Its inverse will be denoted by $\varrho(p)$ for $p>1$.
\end{definition}

\begin{figure}

\includegraphics{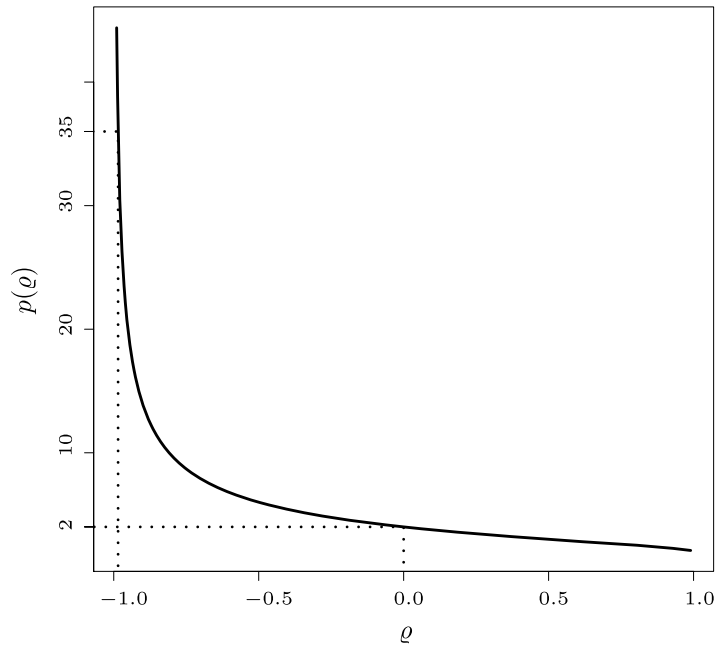}

\caption{The critical curve $p(\varrho), \varrho\in(-1,1)$.}
\label{fig:cl}
\end{figure}

The critical curve is plotted in Figure \ref{fig:cl}. Here, $\varrho
(35)$ and $\varrho(2)$ are marked. Thirty-fifth moments are the key for
%
%
the improved wavespeed result below and the special case $\varrho(2)=0$
is discussed in \cite{AD09}. We will see in Section \ref{sec:moments}
that this curve is closely connected with the exit distribution of
$(B^1_\tau, B^2_\tau)$ from the upper right quadrant which appeared in
Proposition \ref{prop:convlaw} above. The first main theorem states
that the critical curve separates two regimes (independently of $\kappa
$): that of bounded moments and that of unbounded moments.
\begin{theorem}\label{thm:mc}
Suppose $(u_t,v_t)$ is a symbiotic branching process $($either
nonspatial, continuous space or discrete space in arbitrary
dimension$)$ with initial conditions $u_0=v_0=\one$. If $\varrho\in
(-1,1)$, then, for any $\kappa>0$, the following hold for $p>1$:
\begin{longlist}
\item In the recurrent case,
\[
\quad\varrho<\varrho(p)\quad\Leftrightarrow\quad\E
^{1,1}[u_t^p], \E^{\one,\one}[u_t(k)^p]
\mbox{ and }
\E^{\one
,\one}
[u_t(x)^p ] \mbox{ are bounded in }t.
\]
\item In the transient case,
\[
\varrho<\varrho(p) \quad\Rightarrow\quad\E^{\one,\one}
[u_t(k)^p ] \mbox{ is bounded in }t.
\]
\end{longlist}
Due to symmetry the same holds for $\E^{1,1}[v_t^p]$, $\E^{\one,\one
}[v_t(k)^p]$ and $\E^{\one,\one}[v_t(x)^p]$.
\end{theorem}

Note that the theorem provides information about all positive real
moments, not just integer moments. In the area below the critical curve
in Figure \ref{fig:cl}, the moments remain bounded. On and above the
critical curve, in the recurrent case, the moments grow to infinity.
\begin{remark}\label{shift}
For $\varrho=-1$ the curve could be extended with $p(-1)=\infty$. In
terms of the previous theorem this makes sense since for $\varrho=-1$,
symbiotic branching processes with initial conditions $u_0=v_0=\one$ are
bounded by $2$. This is justified by a simple observation: for initial
conditions $u_0=v_0\equiv1/2$ symbiotic branching processes with
$\varrho=-1$ are solutions of the stepping stone model and, hence,
bounded by $1$. Uniqueness in law of solutions implies that solutions
$(u_t,v_t)$ with initial conditions $(cu_0,cv_0)$ are equal in law to
solutions $c$ times solutions with initial conditions $(u_0,v_0)$.
\end{remark}

With this first understanding of the effect of $\varrho$ on moments, we
may discuss integer moments for the discrete-space model in more
detail. Let us first recall some known results for solutions $(w_t)$ of
the parabolic Anderson model (see Example \ref{ex2}) where only the
parameter $\kappa$ appears. Using It\^o's lemma, one sees that
$m(t,k_1,\ldots,k_n):=\E^{\one}[w_t(k_1)\cdots w_t(k_n)]$ solves the
(discrete-space) partial differential equation
\[
\frac{\partial}{\partial t}m(t,k_1,\ldots,k_n)=\Delta m(t,k_1,\ldots
,k_n)+V(k_1,\ldots,k_n)m(t,k_1,\ldots,k_n)
\]
with homogeneous initial conditions. Here, the potential $V$ is given
by
\[
V(k_1,\ldots,k_n)=\kappa\sum_{1\leq i<j\leq n}\delta_0(k_i-k_j).
\]
Since $H=-\Delta-V$ is an $n$-particle Schr\"odinger operator, many
properties are known from the physics literature. In particular, it is
well known that in the recurrent case (the potential is nonnegative)
exponential growth of solutions holds for any $\kappa>0$. By contrast,
in the transient case the discrete Laplacian requires a stronger
perturbation before we see exponential growth. Intuitively from the
particle picture this should be true since the potential $V$ only
increases solutions if particles meet, which occurs less frequently in
the transient case. For the transient case (see, e.g.,
\cite{CM94} or \cite{GdH07} for more precise results), there is a decreasing
sequence $\kappa(n)$ such that
\[
\E^{\one}[w_t(k)^n]\mbox{ is bounded in }t
\quad\Leftrightarrow\quad\kappa<\kappa(n)
\]
and for the Lyapunov exponents
\[
\gamma_n(\kappa):=\lim_{t\rightarrow\infty}\frac{1}{t}\log\E
^{\one
}[w_t(k)^n]>0 \quad\Leftrightarrow\quad\kappa>\kappa(n).
\]
These results can be proved with the $n$-particle path-integral
representation in which solutions are expressed as
\[
m(t,k_1,\ldots,k_n)=\E\bigl[e^{\kappa\int_0^t V(X_s^1,\ldots,X_s^n) \,ds} \bigr],
\]
where $(X^1_t),\ldots,(X^n_t)$ are independent simple random walks started
in\break $k_1,\ldots,k_n$.

Coming back to the symbiotic branching model, we ask whether or not the
$n$th Lyapunov exponents
\[
\gamma_{n}(\varrho,\kappa):=\lim_{t\rightarrow\infty}\frac
{1}{t}\log\E
^{\one,\one}[u_t(k)^{n}]
\]
exist and in which cases $\gamma_n(\varrho,\kappa)$ is strictly
positive. As for the parabolic Anderson model, there is a system of
partial differential equations describing the moments (see Proposition
16 of \cite{EF04} for the continuous-space model) and an $n$-particle
path-integral representation of the moments. In addition to the
independent motion, the particles now carry a color which randomly
changes if particles of the same color stay at the same site (see Lemma
\ref{la:mdual}). With $L_t^=$ denoting collision times of particles of
same color and $L_t^{\neq}$ denoting collision times of particles of
different colors, the path-integral representation of moments reads
\[
\E^{\one,\one}[u_t(k)^n]=\E\bigl[e^{\kappa(L_t^=+\varrho L_t^{\neq})} \bigr].
\]
This representation is more involved than the path-integral
representation for the parabolic Anderson model since, in addition to
the motion of particles, a second stochastic mechanism is included.
Nonetheless, we use it to prove the following theorem which reveals
that even in the recurrent case a nontrivial transition occurs.
\begin{theorem}\label{thm:im}
For solutions of $\mathrm{dSBM}(\varrho,\kappa)_{\one,\one}$, in any
dimension, the following hold for $n\in\N, n>1$:
\begin{longlist}
\item$\gamma_{n}(\varrho,\kappa)$ exists for any $\varrho
\in[-1,1]$, $\kappa>0$,
\item$\gamma_{n}(\varrho(n),\kappa)=0$ for any $\kappa>0$,
\item for any $\varrho>\varrho(n)$ there is a critical
$\kappa(n)$ such that $\gamma_{n}(\varrho,\kappa)>0$ if $\kappa
>\kappa(n)$.
\end{longlist}
\end{theorem}

Combined with Theorem \ref{thm:mc}, parts (ii) and (iii) emphasize the
``criticality'' of the critical curve: for $\varrho<\varrho(n)$,
moments stay bounded, for $\varrho=\varrho(n)$ moments grow
subexponentially fast to infinity, and for $\varrho>\varrho(n)$ moments
grow exponentially fast if $\kappa$ is large enough.
\begin{remark}
As discussed above, for the parabolic Anderson model it is natural that
in the transient case perturbing the critical case does not immediately
yield exponential growth, whereas perturbing the recurrent case does
immediately lead to exponential growth. It is clear that in the
transient case the gap in (iii) of Theorem \ref{thm:im} is really
necessary: for small $\kappa$ moments of the parabolic Anderson model
are bounded. Since moments of symbiotic branching are dominated by
moments of the parabolic Anderson model (see Lemma \ref{la:mdual}), for
small $\kappa$ moments are bounded for all $\varrho$.
\end{remark}

In the case $p\notin\N$ there seems to be no reason why exponential
growth should fail. Unfortunately, in this case there is no moment
duality and hence the most useful tool to analyze exponential growth is
not available.
\begin{conj}
\label{conj:1}
In the recurrent case the moment diagram for symbiotic branching
(Figure \ref{fig:cl}) describes the moments as follows: pairs
$(\varrho
,p)$ below the critical curve correspond precisely to bounded moments,
pairs at the critical curve correspond to moments which grow
subexponentially fast to infinity and pairs above to the critical curve
correspond to exponentially growing moments.
\end{conj}

A deeper understanding of the Lyapunov exponents as functions of
$\varrho,\kappa$ remains mainly open (for an upper bound see
Proposition \ref{up}). For second moments [$\varrho(2)=0$] this is
carried out in \cite{AD09}. It is shown that exponential growth holds
for $\varrho>0$ and arbitrary $\kappa>0$ in the recurrent case, whereas
only for $\kappa>2/(\varrho G_{\infty}(0,0))$ in the transient case.
Here $G_{\infty}$ denotes the Green function of the simple random walk.
The exponential (and subexponential) growth rates were analyzed in
detail by Tauberian theorems.

A direct application of Theorem \ref{thm:im} is so-called intermittency
of solutions. One says a spatial system with Lyapunov exponents $\gamma
_p$ is $p$-intermittent if
\[
\frac{\gamma_p}{p}<\frac{\gamma_{p+1}}{p+1}.
\]
Intermittent systems concentrate on few peaks with extremely high
intensity (see~\cite{GM90}). The results above show that as $\varrho$
tends to $-1$, solutions (at least for large $\kappa$) are
$p$-intermittent for $p$ tending to infinity. This holds since for
fixed $\varrho$, the $p$th moments are bounded if $(\varrho,p)$ lies
below the critical curve. Increasing $p$ (and $\kappa$ if necessary)
there is a first $p$ such that the $p$th Lyapunov exponent is positive.
Intermittency for higher exponents suggests that the effect gets
weaker. This is to be expected since for $\varrho=-1$ solutions with
homogeneous initial conditions are bounded and, hence, solutions do not
produce high peaks at all. Making this effect more precise, in
particular combined with the effect of Proposition \ref{prop:convlaw},
is an interesting task for the future.

\subsection{Speed of propagation of the interface} Let us conclude with
a direct application of the moment bounds. Here, we will be concerned
with an improved upper bound on the speed of the propagation of the
interface of continuous-space symbiotic branching processes which
served to some extent as the motivation for this work. To explain this,
we need to introduce the notion of the interface of continuous-space
symbiotic branching processes introduced in \cite{EF04}.
\begin{definition}\label{def:ifc}
The interface at time $t$ of a solution $(u_t,v_t)$ of the symbiotic
branching model $\mathrm{cSBM}(\varrho,\kappa)_{u_0, v_0}$ with
$\varrho\in[-1,1]$ is defined as
\[
\mathrm{Ifc}_t = \operatorname{cl} \{x\dvtx u_t(x) v_t(x) > 0 \},
\]
where $\operatorname{cl}(A)$ denotes the closure of the set $A$ in $\R$.
\end{definition}

In particular, we will be interested in complementary Heaviside
initial conditions
\[
u_0(x) = \mathbf{1}_{\R^-}(x) \quad\mbox{and}\quad v_0(x) = \mathbf{1}_{\R
^+}(x),\qquad x \in\R.
\]
The main question addressed in \cite{EF04} is whether for the above
initial conditions the so-called compact interface property holds, that
is, whether the interface is compact at each time almost surely. This
is answered affirmatively in Theorem 6 in \cite{EF04}, together with
the assertion that the interface propagates with at most linear speed,
that is, for each $\varrho\in[-1,1]$ there exists a constant $c>0$
and a finite random-time $T_0$ so that almost surely for all $T \ge T_0$
\[
\bigcup_{t \le T} \operatorname{Ifc}_t \subseteq[-cT, cT ].
\]
Heuristically, due to the scaling property of the symbiotic branching
model (Lemma 8 of \cite{EF04}) one expects that the interface should
move with a square-root speed. Indeed, with the help of Theorem \ref
{thm:mc} one can strengthen their result, at least for sufficiently
small $\varrho$,
to obtain almost square-root speed.
\begin{theorem}\label{cor:wavespeed}
Suppose $(u_t,v_t)$ is a solution of $\mathrm{cSBM}(\varrho,\kappa
)_{1_{\R^-}, 1_{\R^+}}$ with $\varrho< \varrho(35)$ and $\kappa>0$.
Then there is a constant $C>0$ and a finite random-time $T_0$ such that
almost surely
\[
\bigcup_{t\leq T} \operatorname{Ifc}_t \subseteq\bigl[-C\sqrt{T\log
(T)},C\sqrt{T\log(T)} \bigr]
\]
for all $T>T_0$.
\end{theorem}

The restriction to $\varrho<\varrho(35)$ is probably not necessary and
only caused by the technique of the proof. Though $\varrho(35)\approx
-0.9958$ is rather close to $-1$, the result is interesting. It shows
that sub-linear speed of propagation is not restricted to situations in
which solutions are uniformly bounded as they are for $\varrho=-1$. The
proof is based on the proof of \cite{EF04} for linear speed which
carries over the proof of \cite{T95} for the stepping stone model to
nonbounded processes. We are able to strengthen the result by using a
better moment bound which is needed to circumvent the lack of uniform
boundedness.
%
\begin{remark}\label{r}
We believe that, at least for $\varrho\leq0$, the speed of
propagation should be at most $C'\sqrt t$, for some suitable constant
$C'$, that is, for all $T$ greater than some $T'>0$,
\[
\bigcup_{t\leq T} \operatorname{Ifc}_t \subseteq\bigl[-C'\sqrt
{T},C'\sqrt
{T} \bigr].
\]
However, it seems unclear how to obtain such a refinement of
Theroem \ref{cor:wavespeed} based on our moment results and the
method of \cite{T95} (resp., \cite{EF04}). As subexponential bounds of
higher moments
cannot be avoided (see the proof of the fluctuation term estimate
Lemma \ref{la:ma}),
our results on the behavior of higher moments show that at present, in
light of Conjecture \ref{conj:1},
one can only hope for stronger results for very small $\varrho$.

To overcome this limitation, new methods need to be employed. The
authors think that a possible approach
could be based on the scaling property (Lemma 8 of \cite{EF04}) and
recent results by Klenke and Oeler \cite{KO09}.
Recall that the scaling property states that if $(u_t, v_t)$ is a
solution to $\operatorname{cSBM}(\varrho, \kappa)_{u_0, v_0}$,
then
\[
(u_t(x)^K, v_t(x)^K ) := \bigl(u_{Kt}\bigl(\sqrt K x\bigr), v_{Kt}\bigl(\sqrt K
x\bigr) \bigr),\qquad x \in\R, K >0,
\]
is a solution to $\operatorname{cSBM}(\varrho, K \cdot\kappa
)_{u^K_0, v^K_0}$
(with suitably transformed initial states $u^K_0, v^K_0$).
In other words, a diffusive time--space rescaling leads to the
original model with a suitably increased branching rate $\kappa$.
Klenke and Oeler \cite{KO09} show that, at least for the mutually
catalytic model
in discrete space, a nontrivial limiting process as $\kappa\to\infty
$ exists.
This limit is called ``infinite rate mutually catalytic branching process''
(see also \cite{KM09a,KM09b} for a further discusion). In
particular, in
Corollary 1.2 of \cite{KO09} they claim that, under suitable
assumptions, a nontrivial interface for the limiting process exists,
which would in turn predict a square-root speed of propagation in our case.
However, to make this approach rigorous is beyond the scope of the
present paper.
\end{remark}
\begin{remark}[(Shape of the interface)]
Note that our results give only limited information about the shape of
the interface.
For the case $\varrho=-1$, that is, with locally constant total
population size,
it is shown in \cite{MT97} that there exists a unique stationary
interface law, which may therefore be interpreted
as a ``stationary wave'' whose position fluctuates at the boundaries,
according to \cite{T95}, like a Brownian motion, hence explaining the
square-root speed (note that for both results, suitable bounds on
fourth mixed moments are required).
However, for $\varrho> -1$, the population sizes of the interface are
expected to fluctuate significantly and it seems unclear how this affects
the shape and speed of the interface, in particular the formation of a
``stationary wave.''
The significance of fourth mixed moments
might even lead to a phase-transition in $\varrho$.
This gives rise to many interesting open questions.
\end{remark}

\section{Basic properties and duality}\label{sec:bpd}
In this section we review the setting and properties of the
discrete-space model, whereas for continuous-space we refer
to~\cite{EF04}. Note that instead of using the state space of tempered
functions alternatively we may use a suitable Liggett--Spitzer space.
As the results are only presented for the discrete Laplacian this does
not play a crucial role. For a discussion of the mutually catalytic
branching model in the Liggett--Spitzer space see \cite{CDG04}.

\subsection{Basic properties}
\label{subsec:bp}
For functions $f,g\dvtx\Z^d \to\R$ we abbreviate $\langle
f,g\rangle
=\sum
_kf(k)g(k)$. With $\phi_{\lambda}(k)=e^{\lambda|k|}$ the space of pairs
of tempered sequences is defined by
\[
M^2_{\mathrm{tem}}= \{(u,v) | u,v\dvtx\Z^d\rightarrow\R_{\geq0} ,
\langle u,\phi_{\lambda}\rangle,\langle v,\phi_{\lambda}\rangle
<\infty\
\forall\lambda<0 \}.
\]
The space of continuous paths is denoted by
\[
\Omega_{\mathrm{tem}}=C(\R_{\geq0},M_{\mathrm{tem}}^2).
\]
Similarly, the space of pairs of rapidly decreasing sequences is
defined by
\[
M^2_{\mathrm{rap}}= \{(u,v) | u,v\dvtx\Z^d\rightarrow\R_{\geq0} ,
\langle u,\phi_{\lambda}\rangle,\langle v,\phi_{\lambda}\rangle
<\infty\
\forall\lambda>0 \}
\]
and the corresponding path space by
\[
\Omega_{\mathrm{rap}}=C(\R_{\geq0},M_{\mathrm{rap}}^2).
\]
Weak solutions are defined as in \cite{DP98} for $\varrho=0$. In much
the same way as for Theorems 1.1 and 2.2 of \cite{DP98}, we
obtain existence and the Green-function representation.
\begin{proposition}\label{prop:bp1}
Suppose $(u_0,v_0)\in M_{\mathrm{tem}}^2$ (resp., $M^2_{\mathrm{rap}}$), $\varrho\in
[-1,1]$ and $\kappa>0$. Then there is a weak solution of $\mathrm
{dSBM}(\varrho,\kappa)_{u_0,v_0}$ such that $(u_t,v_t)\in\Omega_{\mathrm{tem}}$
(resp., $\Omega_{\mathrm{rap}}$) and for all $(\phi,\psi)\in M^2_{\mathrm{rap}}$
(resp., $M^2_{\mathrm{tem}}$)
%
%
\begin{eqnarray}\label{12}
\langle u_t,\phi\rangle&=&\langle u_0,
P_t\phi\rangle+\sum_{j\in\Z^d}\int_0^t P_{t-s}\phi(j)\sqrt
{\kappa
u_s(j)v_s(j)} \,dB^1_s(j),\\
\label{13}
\langle v_t,\psi\rangle&=&\langle v_0,
P_t\psi\rangle+\sum_{j\in\Z^d}\int_0^t P_{t-s}\psi(j)\sqrt
{\kappa
u_s(j)v_s(j)} \,dB^2_s(j),
\end{eqnarray}
where $P_tf(k)=\sum_{j\in\Z^d}p_t(j,k)f(j)$ is the semigroup associated
to the simple random walk.
In particular, we have
%
%
\begin{eqnarray}
\label{14}
u_t(k)&=&P_tu_0(k)+\sum_{j\in\Z^d}\int_0^tp_{t-s}(j,k)\sqrt{\kappa
u_s(j)v_s(j)} \,dB^1_s(j),\\
\label{15}
v_t(k)&=&P_tv_0(k)+\sum_{j\in\Z^d}\int_0^tp_{t-s}(j,k)\sqrt{\kappa
u_s(j)v_s(j)} \,dB^2_s(j).
\end{eqnarray}
The covariation structure of the Brownian motions is given by
(\ref{eq:dcv}).
\end{proposition}

In fact, (\ref{12}), (\ref{13}) can be seen as the discrete-space
versions of the martingale problem of Definition 3 in \cite{EF04}.
Further, (\ref{14}), (\ref{15}) are the discrete-space versions of the
convolution form given in Corollary 20 of \cite{EF04}.

For the proofs of the longtime behavior of laws and moments, the key
step is to transfer to the total mass processes $\langle u_t,\one
\rangle,
\langle v_t,\one\rangle$. To this end, in a similar way to Proposition
\ref{prop:bp1}, we define
\[
M_{F}^2= \{(u,v) | u,v\dvtx\Z^d\rightarrow\R_{\geq0},
\langle u,1\rangle,\langle v,1\rangle<\infty\}
\]
and
\[
\Omega_{F}=C(\R_{\geq0},M_{F}^2).
\]
For summable initial conditions we obtain the following crucial
martingale characterization.
\begin{proposition}\label{prop:tmmart}
If $(u_0,v_0)\in M^2_F$, then each solution of $\mathrm{dSBM}(\varrho
,\kappa)_{u_0,v_0}$ has the following properties:
$(u_t,v_t)\in\Omega_F$ and $\langle u_t,\one\rangle, \langle
v_t,\one
\rangle$ are nonnegative, continuous, square-integrable martingales
with square-functions
\[
[\langle u_{\cdot},\one\rangle]_t= [\langle v_{\cdot},\one
\rangle]_t=\kappa\int_0^t\langle u_s,v_s\rangle \,ds
\]
and
\[
[\langle u_{\cdot},\one\rangle,\langle v_{\cdot},\one\rangle
]_t=\varrho\kappa\int_0^t\langle u_s,v_s\rangle \,ds.
\]
\end{proposition}

We omit the proofs since they are basically standard. The only step
where one needs to be careful is the existence proof. As usual for such
models one first restricts the space to bounded subsets (boxes) of $\Z
^d$, where standard Markov process theory applies. Enlarging the boxes
one obtains a sequence of processes which are shown to converge to a
limiting process solving $\mathrm{dSBM}(\varrho,\kappa)_{u_0,v_0}$. To
prove tightness of the approximating sequence, the moments need to be
bounded uniformly in the size of the boxes. Here, more care than for
$\varrho=0$ in \cite{DP98} is needed. The uniform moment bound can, for
instance, be achieved using a colored particle moment duality for each
box similar to the one of Lemma \ref{la:mdual}.

\subsection{Dualities}\label{subsec:dualities}

The symbiotic branching model exhibits an exceptionally rich duality
structure, providing powerful tools for the analysis of the longtime properties.

\subsubsection{Colored particle moment dual}\label{subsec:cpd}

We now recall the two-colors particle moment-duality introduced in
Section 3.1 of \cite{EF04}. Since the dual Markov process is presented
rigorously in \cite{EF04} we only sketch the pathwise behavior. To find
a suitable description of the mixed moment
\[
\E^{u_0,v_0}[u_t(k_1)\cdots u_t(k_{n})v_t(k_{n+1})\cdots v_t(k_{n+m})],
\]
$n+m$ particles are located in $\Z^d$. Each particle moves as a
continuous-time simple random walk independent of all other particles.
At time $0$, $n$ particles of color $1$ are located at positions
$k_1,\ldots,k_n$ and $m$ particles of color $2$ are located at positions
$k_{n+1},\ldots,k_{n+m}$. For each pair of particles, one of the pair
changes color when the time the two particles have spent in the same
site, while both have same color, first exceeds an (independent)
exponential time with parameter $\kappa$. Let
\begin{eqnarray*}
L_t^=&=&\mbox{total collision time of all pairs of same colors up to
time }t,\\
L_t^{\neq}&=&\mbox{total collision time of all pairs of different
colors up to time }t,\\
l^1_t(a)&=&\mbox{number of particles of color }1\mbox{ at site }a\mbox{
at time $t$},\\
l^2_t(a)&=&\mbox{number of particles of color }2\mbox{ at site }a\mbox{
at time $t$},\\
(u_0,v_0)^{l_t}&=&\prod_{a\in\Z^d}u_0(a)^{l_t^1(a)}v_0(a)^{l_t^2(a)}.
\end{eqnarray*}
Note that since there are only $n+m$ particles, the infinite product is
actually a finite product and hence well defined. The following lemma
is taken from Section~3 of~\cite{EF04}.
\begin{lemma}\label{la:mdual}
Let $(u_t,v_t)$ be a solution of $\mathrm{dSBM}(\varrho,\kappa
)_{u_0,v_0}$, $\kappa>0$ and $\varrho\in[-1,1]$. Then, for any
$k_i\in\Z
^d$, $t\geq0$,
\[
\E^{u_0,v_0}[u_t(k_1)\cdots u_t(k_{n})v_t(k_{n+1})\cdots
v_t(k_{n+m})]=\E\bigl[(u_0,v_0)^{l_t}e^{\kappa(L_t^=+\varrho L_t^{\neq
})} \bigr],
\]
where the dual process behaves as explained above.
\end{lemma}

Note that for homogeneous initial conditions $u_0=v_0=\one$, the first
factor in the expectation of the right-hand side equals $1$. In the
special case $\varrho=1$, \mbox{$u_0=v_0=\one$} Lemma \ref{la:mdual} was already
stated in \cite{CM94}, reproved in \cite{GdH07} and used to analyze the
Lyapunov exponents of the parabolic Anderson model.

For $\varrho\neq1$, the difficulty of the dual process is based on the
two stochastic effects: on the one hand, one has to deal with collision
times of random walks which were analyzed in \cite{GdH07};
additionally, particles have colors either $1$ or $2$ which change dynamically.
\begin{remark}\label{rdual}
Similar dualities hold for $\mathrm{cSBM}$ and $\mathrm{SBM}$. For
continuous-space, the random walks are replaced by Brownian motions and
the collision times of the random walks by collision local times of the
Brownian motions (see Section~4.1 in \cite{EF04}). The simplest case is
the nonspatial symbiotic branching model where the particles stay at
the same site and local times are replaced by real times (see Theorem
3.2 of \cite{R95} or Proposition A5 of \cite{DFX05}).
\end{remark}

\subsubsection{Self-duality}\label{subsec:sduality}

Mytnik \cite{M99} introduced a self-duality for the continuous-space
mutually catalytic branching model to obtain uniqueness of solutions of
the corresponding martingale problem. This can be extended to symbiotic
branching models for $\varrho\in(-1,1)$ as shown in Proposition 5 of
\cite{EF04}. The discrete-space self-duality for $\varrho=0$ was proved
in Theorem 2.4 of \cite{DP98}. We first need more spaces of sequences:
\[
E= \{ (x,y)\dvtx(x,|y|)\in M^2_{\mathrm{tem}}, |y(k)|\leq x(k)\ \forall k\in\Z
^d \}
\]
and
\[
\tilde E= \{ (x,y)\in E\dvtx x\in M_{\mathrm{rap}} \}\supset\{ (x,y)\in
E\dvtx x\mbox{ has bounded support} \}=\tilde{E}_f.
\]
In the sequel, the space $E$ and its subspaces will be used for
$(x,y)=(u_t+v_t,u_t-v_t)$. The duality function for $\varrho\in(-1,1)$
maps $ E \times\tilde E$ to $\mathbb{C}$ via
%
%
\begin{equation}\label{eq:selfdf2}
H (u,v,\tilde{u}, \tilde{v} )=\exp\bigl(-\sqrt{1-\varrho}\langle
u,\tilde{u} \rangle+i\sqrt{1+\varrho}\langle v,\tilde{v} \rangle\bigr).
\end{equation}
With this definition the generalized Mytnik duality states:
\begin{lemma}\label{la:sduality}
For $\varrho\in(-1,1)$, $\kappa>0$, $(u_0,v_0)\in M^2_{\mathrm{tem}}$ and
$(\tilde u_0, \tilde v_0)\in M^2_{\mathrm{rap}}$ let $(u_t,v_t)$ be a solution
of $\mathrm{dSBM}(\varrho,\kappa)_{u_0, v_0}$ and $(\tilde
{u}_t,\tilde
{v}_t)$ be a solution of $\mathrm{dSBM}(\varrho,\kappa)_{\tilde u_0,
\tilde v_0}$. Then the following holds:
\begin{eqnarray*}
&&
\E^{u_0,v_0} [H(u_t+v_t,u_t-v_t,\tilde u_0+\tilde v_0, \tilde
u_0-\tilde v_0) ]\\
&&\qquad=\E^{\tilde{u}_0,\tilde{v}_0} [H(u_0+v_0,u_0-v_0,\tilde{u}_t+\tilde
v_t,\tilde u_t-\tilde{v}_t) ].
\end{eqnarray*}
\end{lemma}

Analogously, the self-duality relation holds for the nonspatial model
with duality function
\[
H^0(u,v,\tilde u,\tilde v)=\exp\bigl(-\sqrt{1-\varrho}u\tilde{u}+i\sqrt
{1+\varrho}v\tilde{v} \bigr),
\]
mapping $(\R_{\geq0}\times\R_{\geq0})^2$ to $\mathbb{C}$.


\section{Weak longtime convergence}\label{sec:comnvlaw}
In this section we discuss weak longtime convergence of symbiotic
branching models and prove Proposition \ref{prop:convlaw}. We proceed
in two steps: first, we prove convergence in law to some limit law
following the proof of \cite{DP98} for $\varrho=0$. Second, to
characterize the limit law for the spatial models, we reduce the
problem to the nonspatial model.
\begin{proposition}\label{prop:wconv}
Let $\varrho\in(-1,1), \kappa>0$ and $(u_t,v_t)$ a solution of
either $\mathrm{cSBM}(\varrho,\kappa)_{u_0, v_0}$ or $\mathrm
{dSBM}(\varrho,\kappa)_{u_0, v_0}$
with initial conditions $u_0=\bu, v_0=\bv$. Then, as $t \to\infty$,
the law of $(u_t,v_t)$ converges weakly on $M_{\mathrm{tem}}^2$ to some limit
$(u_{\infty},v_{\infty})$.
\end{proposition}
\begin{pf}
The proof is only given for the discrete spatial case and the
continuous case is completely analogous. Let us first recall the
strategy of \cite{DP98} for $\varrho=0$ which can also be applied with
the generalized self-duality required here. Convergence of $(u_t,v_t)$
in $M_{\mathrm{tem}}^2$ follows from convergence of $(u_t+v_t,u_t-v_t)$ in $E$.
Using Lemma 2.3(c) of \cite{DP98}, it suffices to show convergence of
$\E^{\bu,\bv}[H(u_t+v_t,u_t-v_t,\phi, \psi)]$ for all $(\phi,\psi
)\in
\tilde E_f$. Furthermore, the limit $(u_{\infty},v_{\infty})$ is
uniquely determined by $\E^{\bu,\bv}[H(u_{\infty}+v_{\infty
},u_{\infty
}-v_{\infty},\phi,\psi)]$ (see Lemma 2.3(b) of \cite{DP98}). Hence, it
suffices to show convergence of
%
%
\begin{equation}\label{456}
\E^{\bu,\bv}[H(u_t+v_t,u_t-v_t,\phi,\psi)]=\E^{\bu,\bv}
\bigl[e^{-\sqrt
{1-\varrho}\langle u_t+v_t,\phi\rangle+i\sqrt{1+\varrho}\langle
u_t-v_t,\psi\rangle} \bigr],\hspace*{-32pt}
\end{equation}
for all $(\phi,\psi)\in\tilde{E}_f$.
Note that the technical condition of Lemma 2.3(c) of \cite{DP98} is
fullfilled since due to Proposition \ref{prop:bp1}
\[
\E^{\bu,\bv}[\langle u_t+v_t,\phi_{-\lambda}\rangle]=(u+v)\langle
\one
,P_t \phi_{-\lambda}\rangle<C<\infty.
\]
To ensure convergence of (\ref{456}) we employ the generalized Mytnik
self-duality of Lemma \ref{la:sduality} with $\tilde{u}_0:=\frac
{\phi
+\psi}{2},\tilde{v}_0:=\frac{\phi-\psi}{2}$:
%
%
\begin{eqnarray}\label{eq:dualconv}
&&\E^{\bu,\bv} \bigl[e^{-\sqrt{1-\varrho} \langle u_t+v_t,\phi\rangle
+i\sqrt{1+\varrho}
\langle u_t-v_t,\psi\rangle} \bigr] \nonumber\\
&&\qquad=\E^{u_0,v_0} \bigl[e^{-\sqrt{1-\varrho}\langle u_t+v_t,\tilde
{u}_0+\tilde{v}_0\rangle
+i\sqrt{1+\varrho}\langle u_t-v_t,\tilde{u}_0-\tilde{v}_0\rangle} \bigr]
\nonumber\\[-8pt]\\[-8pt]
&&\qquad=\E^{\tilde{u}_0,\tilde{v}_0} \bigl[e^{-\sqrt{1-\varrho}\langle u_0+v_0,
\tilde{u}_t+\tilde{v}_t\rangle+i\sqrt{1+\varrho}\langle u_0-v_0,
\tilde{u}_t-\tilde{v}_t\rangle} \bigr]\nonumber\\
&&\qquad=\E^{\tilde{u}_0,\tilde{v}_0} \bigl[e^{-\sqrt{1-\varrho}(u+v)\langle
\mathbf{1},
\tilde{u}_t+\tilde{v}_t\rangle+i\sqrt{1+\varrho}(u-v)\langle
\mathbf{1},
\tilde{u}_t-\tilde{v}_t\rangle} \bigr].\nonumber
\end{eqnarray}
By assumption, $\tilde u_0, \tilde v_0$ have compact support and hence
by Proposition \ref{prop:tmmart} the total-mass processes
$\langle\mathbf{1}, \tilde{u}_t \rangle$ and $\langle\mathbf{1},
\tilde
{v}_t \rangle$ are
nonnegative martingales. By the martingale convergence theorem
$\langle\mathbf{1}, \tilde{u}_t \rangle$ and $\langle\mathbf{1},
\tilde
{v}_t \rangle$ converge almost surely to finite limits denoted by
$\langle\mathbf{1}, \tilde{u}_{\infty} \rangle$, $\langle\mathbf
{1},\tilde
{v}_{\infty} \rangle$. Finally, the dominated convergence theorem
implies convergence of the right-hand side of (\ref{eq:dualconv}) to
%
%
\begin{equation}\label{eq:elimit}
\E^{\tilde u_0,\tilde v_0} \bigl[e^{-\sqrt{1-\varrho}(u+v)\langle
\mathbf{1},\tilde{u}_{\infty}+\tilde{v}_{\infty}\rangle+i\sqrt
{1+\varrho}(u-v)
\langle\mathbf{1},\tilde{u}_{\infty}-\tilde{v}_{\infty}\rangle} \bigr].
\end{equation}
Combining the above, we have proved convergence of
\[
\E^{\bu,\bv} \bigl[e^{-\sqrt{1-\varrho} \langle u_t+v_t,\phi\rangle
+i\sqrt{1+\varrho}\langle u_t-v_t,\psi\rangle} \bigr],
\]
which ensures weak convergence of $(u_t,v_t)$ in $M_{\mathrm{tem}}^2$ to some
limit which is uniquely determined by (\ref{eq:elimit}).
\end{pf}

Again, as in Remark \ref{all}, the previous proposition can be proved
for nondeterministic initial conditions as in \cite{CKP00}.

The rest of this section is devoted to identifying the limit
$(u_{\infty
},v_{\infty})$ in the recurrent case. Before completing the proof of
Theorem \ref{prop:convlaw} we discuss a version of Knight's extension
of the Dubins--Schwarz theorem (see \cite{KS98}, 3.4.16) for
nonorthogonal continuous local martingales.
\begin{lemma}\label{la:eds}
Let $(N_t)$ and $(M_t)$ be continuous local martingales with
$N_0=M_0=0$ almost surely. Assume further that, for $t \ge0$,
\[
[ M_{\cdot},M_{\cdot}]_t=[ N_{\cdot},N_{\cdot}]_t \quad\mbox{and}\quad
[ M_{\cdot},N_{\cdot}]_t=\varrho[ M_{\cdot},M_{\cdot}]_t\qquad
\mbox{a.s.},
\]
where $\varrho\in[-1,1]$. If $[ M_{\cdot},M_{\cdot}]_{\infty
}=\infty$
a.s., then
\[
(B^1_t,B^2_t):=\bigl(M_{T(t)},N_{T(t)}\bigr)
\]
is a pair of Brownian motions with covariances $[ B^1_{\cdot
},B^2_{\cdot
}]_t=\varrho t$, where
%
%
\begin{equation}\label{eq:timeshift}
T(t)=\inf\{s\dvtx[ M_{\cdot},M_{\cdot}]_s>t \}.
\end{equation}
\end{lemma}
\begin{pf}
It follows from the Dubins--Schwarz theorem that $B^1,B^2$ are each
Brownian motions. Further, by the definition of $T(t)$ we obtain the claim
\[
[ B^1_{\cdot}, B^2_{\cdot}]_t=[ M_{\cdot},
N_{\cdot}]_{T(t)}=\varrho[ M_{\cdot},M_{\cdot}]_{T(t)}=\varrho t.
\]
\upqed\end{pf}
\begin{remark}
If $T^*:=[ M_{\cdot},M_{\cdot}]_{\infty}<\infty$ the situation becomes
slightly more delicate but one can use a local version of Lemma \ref
{la:eds}. Indeed, define, for $t \ge0$,
%
%
\begin{equation}\label{eq:delicate}
B^1_t:=
\cases{
M_{T(t)}, &\quad for $t < T^*$,\cr
M_{T^*}, &\quad for $t \geq T^*$,}
\end{equation}
where the time-change $T$ is given in (\ref{eq:timeshift}) and define
$B^2$ analogously for $N$ (recall that $[M_{\cdot},M_{\cdot}]_t=[
N_{\cdot},N_{\cdot}]_t$). Then the processes $B^1, B^2$ are Brownian
motions stopped at time $T^*$.
The covariance is again given by
\[
[ B^1_{\cdot}, B^2_{\cdot}]_{t\wedge T^*} = \varrho(t\wedge T^*),\qquad
t \ge0.
\]
\end{remark}

For the rest of this section let $B^1,B^2$ be standard Brownian
motions with covariance
%
%
\begin{equation}
[ B^1_{\cdot},B^2_{\cdot}]_t=\varrho t
\end{equation}
started in $u,v$, denote their expectations by $E^{u,v}$, and let
\[
\tau=\inf\{t\dvtx B^1_tB^2_t=0 \}.
\]

The above discussion can now be used to understand the longtime
behavior of symbiotic branching processes. We start by giving a proof
for the nonspatial symbiotic branching model and then modify the proof
to capture the corresponding result for the spatial models.
\begin{proposition}\label{prop:sconv}
Let $(u_t,v_t)$ be a solution of $\mathrm{SBM}(\varrho,\kappa)_{u,v}$.
Then, as $t\to\infty$, $(u_t,v_t)$ converges almost surely to some
$(u_{\infty},v_{\infty})$. Furthermore, $\mathcal{L}^{u,v}(u_{\infty
},v_{\infty})=\mathcal{L}^{u,v}(B^1_{\tau},B^2_{\tau})$ with
$B^1_{\tau
},B^2_{\tau}$ from Proposition \ref{prop:convlaw}.
\end{proposition}
\begin{pf}
Solutions of the nonspatial symbiotic branching model are nonnegative
martingales and hence converge almost surely. This implies the first
part of the claim and it only remains to characterize the limit.
Obviously, the $L^2$-martingales $(u_t),(v_t)$ satisfy the
cross-variation structure assumptions of Lemma~\ref{la:eds} and, thus,
$(u_{t},v_{t})= (B^1_{T^{-1}(t)}, B^2_{T^{-1}(t)} )$. To obtain
the result, we need to check that $T^{-1}(\infty)= \tau$. By definition
of $\mathrm{SBM}$, the time-change is given by
%
%
\begin{equation}\label{eq:invtime}
T^{-1}(t) = [ u_\cdot, u_\cdot]_t = \biggl[ \int_0^\cdot\sqrt{\kappa
u_s u_s} \,dB^1_s, \int_0^\cdot\sqrt{\kappa u_s v_s} \,dB^1_s \biggr]_t
= \kappa\int_0^t u_s v_s \,ds.\hspace*{-32pt}
\end{equation}
To see that $T^{-1}(\infty)=\tau<\infty$, first note that $T^{-1}(t)
\le\tau$ for all $t\geq0$. This is true since $u_t=B^1_{T^{-1}(t)},
v_t=B^2_{T^{-1}(t)}$ and solutions of $\mathrm{SBM}$ are nonnegative.
To argue that $T^{-1}(t)$ increases to $\tau$, more care is needed.
Since the martingales converge almost surely, $T^{-1}(t)$ converges to
some value $a\leq\tau$. Suppose $a<\tau$, then $(u_t,v_t)$ converges
to some $(x,y)$ with $x,y>0$. This yields a contradiction since
$T^{-1}(t)=\kappa\int_0^tu_sv_s \,ds$ would increase to infinity.
Hence, almost surely,
\[
(u_t,v_t)= \bigl(B^1_{T^{-1}(t)},B^2_{T^{-1}(t)} \bigr) \stackrel
{t\rightarrow
\infty}{\rightarrow}\bigl(B^1_{T^{-1}(\infty)},B^2_{T^{-1}(\infty
)}\bigr)=(B^1_{\tau},B^2_{\tau}).
\]
\upqed\end{pf}

In particular, the proof of Proposition \ref{prop:sconv} provides an
important relation for $(B^1_{\tau},B^2_{\tau})$. As remarked below
Lemma \ref{la:sduality}, the self-duality also works in the nonspatial model:
\begin{eqnarray*}
&&\E^{u_0,v_0}\bigl[e^{-\sqrt{1-\varrho}(u_t+v_t)(\tilde{u}_0+\tilde{v}_0
)+i\sqrt{1+\varrho}(u_t-v_t)(\tilde{u}_0-\tilde{v}_0 )}\bigr]\\
&&\qquad=\E^{\tilde{u}_0,\tilde{v}_0}\bigl[e^{-\sqrt{1-\varrho
}(u_0+v_0)(\tilde
{u}_t+\tilde{v}_t)+i\sqrt{1+\varrho}(u_0-v_0)(\tilde{u}_t-\tilde
{v}_t )}\bigr],
\end{eqnarray*}
where both $(u_t,v_t)$ and $(\tilde{u}_t,\tilde{v}_t)$ are solutions of
$\mathrm{SBM}(\varrho,\kappa)$ with different initial conditions. As
shown in the proof of Proposition \ref{prop:sconv}, $(u_t,v_t)$ [resp.,
$(\tilde{u}_t,\tilde{v}_t)$] converges almost surely to $(B^1_{\tau
},B^2_{\tau})$ with initial condition $(u_0,v_0)$ [resp., $(\tilde
u_0,\tilde v_0$)]. Using dominated convergence, this shows the
following duality relation for $(B^1_{\tau},B^2_{\tau})$ when started
in initial conditions $(u,v)$, $(\tilde{u},\tilde{v})$:
%
%
\begin{eqnarray}\label{dual0}
&&E^{u,v} [H^0(B^1_{\tau}+B^2_{\tau},B^1_{\tau}-B^2_{\tau},\tilde
u+\tilde v,\tilde u-\tilde v) ]\nonumber\\[-8pt]\\[-8pt]
&&\qquad=E^{\tilde{u},\tilde{v}} [H^0(B^1_{\tau}+B^2_{\tau},B^1_{\tau
}-B^2_{\tau},u+v,u-v) ].\nonumber
\end{eqnarray}
\begin{pf*}{Proof of Proposition \ref{prop:convlaw}}
Again, the proof is only presented in the discrete spatial setting
since the continuous case is analogous. We retain the notation of the
proof of Proposition \ref{prop:wconv} where we showed that, as $t$
tends to infinity,
\begin{eqnarray*}
&&\E^{\bu,\bv} \bigl[e^{-\sqrt{1-\varrho} \langle u_t+v_t,\phi\rangle
+i\sqrt{1+\varrho}
\langle u_t-v_t,\psi\rangle} \bigr]\\
&&\qquad\rightarrow\E^{({\phi+\psi})/{2},({\phi-\psi})/{2}}
\bigl[e^{-\sqrt
{1-\varrho}(u+v)\langle\mathbf{1},\tilde{u}_{\infty}+\tilde
{v}_{\infty
}\rangle+i\sqrt{1+\varrho}(u-v) \langle\mathbf{1},
\tilde{u}_{\infty}-\tilde{v}_{\infty}\rangle} \bigr].
\end{eqnarray*}
Let us specify the limit law as for the nonspatial symbiotic branching
process. As seen in Proposition \ref{prop:tmmart} the total-mass
processes $\bar u_t :=\langle\tilde{u}_t,\one\rangle$ and $\bar
v_t:=\langle\tilde{v}_t,\one\rangle$ are nonnegative continuous
$L^2$-martingales with cross-variations $[\bar u_{\cdot}, \bar
v_{\cdot
} ]_t=\varrho[\bar u_{\cdot}, \bar u_{\cdot}]_t=\varrho[\bar
v_{\cdot
}, \bar v_{\cdot}]_t, t\ge0$. Thus, by Lemma \ref{la:eds}, reasoning
as in (\ref{eq:invtime}), $(\bar u_t,\bar v_t)=
(B^1_{T^{-1}(t)},B^2_{T^{-1}(t)} )$, where $B^1,B^2$ are Brownian
motions started in $\bar u_0= \langle\frac{\phi+\psi}{2},\one
\rangle$, $\bar v_0= \langle\frac{\phi-\psi}{2},\one\rangle$ with
covariance $[B^1_{\cdot}, B^2_{\cdot}]_t=\varrho t$ and $T^{-1}(t)
=\kappa\int_0^t \langle u_s,v_s\rangle \,ds$. Again, we need to show
that $T^{-1}(\infty) = \tau$. This is much more subtle than in the
nonspatial case since the quadratic variation might level off even if
both total-mass processes $\bar u_t$, $\bar v_t$ are strictly positive.
In \cite{DP98} it was shown that for $\varrho=0$, almost surely, this
does not happen in the recurrent case [cf. the proof of their
Theorem~1.2(b)]. Their proof can be used directly 
for $\varrho\in(-1,1)$. Hence, almost surely,
%
%
\begin{equation}\label{eq:barconv}
(\langle\tilde{u}_t,\one\rangle,\langle\tilde{v}_t,\one\rangle)
\stackrel
{t\rightarrow\infty}{\rightarrow} (B^1_{\tau}, B^2_{\tau}).
\end{equation}
Combining the above discussion with (\ref{eq:elimit}), we are able to
determine the limit. First, we derived
\begin{eqnarray*}
&&\E^{\bu,\bv} \bigl[e^{-\sqrt{1-\varrho} \langle u_t+v_t,\phi\rangle
+i\sqrt{1+\varrho}\langle u_t-v_t,\psi\rangle} \bigr]\\
&&\qquad\stackrel{t\rightarrow\infty}{\rightarrow}E^{\langle({\phi
+\psi
})/{2},\one\rangle,\langle({\phi-\psi})/{2},\one
\rangle} \bigl[e^{-\sqrt{1-\varrho}(u+v)( B^1_{\tau}+B^2_{\tau})
+i\sqrt
{1+\varrho}(u-v) ( B^1_{\tau}-B^2_{\tau})} \bigr].
\end{eqnarray*}
To use Lemma 2.3(c) of \cite{DP98} we manipulate the right-hand side
using (\ref{dual0}):
\begin{eqnarray*}
&&E^{\langle({\phi+\psi})/{2},\one\rangle,\langle({\phi
-\psi
})/{2},\one
\rangle} \bigl[e^{-\sqrt{1-\varrho}(u+v)( B^1_{\tau}+B^2_{\tau})
+i\sqrt
{1+\varrho}(u-v) ( B^1_{\tau}-B^2_{\tau})} \bigr]\\
&&\qquad=E^{\langle({\phi+\psi})/{2},\one\rangle,\langle({\phi
-\psi
})/{2},\one
\rangle} [H^0 ({B^1_{\tau}}+{B^2_{\tau}},{B^1_{\tau}}-{B^2_{\tau
}},u+v,u-v ) ]\\
&&\qquad=E^{u,v} [H^0 ({B^1_{\tau}}+{B^2_{\tau}},{B^1_{\tau}}-{B^2_{\tau
}},\langle\phi,\one\rangle,\langle\psi,\one\rangle) ]\\
&&\qquad=E^{u,v} [H (\bar{B}^1_{\tau}+\bar{B}^2_{\tau},\bar{B}^1_{\tau
}-\bar{B}^2_{\tau},\phi,\psi) ],
\end{eqnarray*}
where, as in Proposition \ref{prop:convlaw}, $\bar{B}_{\tau}^1$
(resp., $\bar{B}_{\tau}^2$) denotes the constant function
taking only the (random) value $B_{\tau}^1$ (resp., ${B}_{\tau}^2$).
In total we have 
%
\[
\E^{\bu,\bv} [H(u_t+v_t,u_t-v_t,\phi,\psi) ]
\stackrel{t\rightarrow\infty}{\rightarrow}E^{u,v} [H(\bar
{B}^1_{\tau
}+\bar{B}^2_{\tau},\bar{B}^1_{\tau}-\bar{B}^2_{\tau},\phi,\psi) ],
\]
which implies weak convergence in $M_{\mathrm{tem}}^2$ of $(u_t,v_t)$ to $(\bar
{B}^1_{\tau},\bar{B}^2_{\tau})$ by Lemma 2.3(c) of \cite{DP98}.
\end{pf*}

\section{Moments}\label{sec:moments}
In this section we prove Theorems \ref{thm:mc} and \ref{thm:im}. Before
giving the proofs we prove an equivalence for moments of correlated
Brownian motions.

\subsection{Moments of the exit-point and exit-time distribution of
correlated Brownian motions in a quadrant}\label{subsec:exitmoments}
Let $\varrho\in(-1,1)$, $u,v>0$ and $B^1, B^2$ be Brownian motions
started in $u,v$ with
%
%
\begin{equation}\label{eq:bmcor}
\langle B^1_{\cdot}, B^2_{\cdot}\rangle_t=\varrho t.
\end{equation}
The starting points of Brownian motions will be indicated by
superscripts in probabilities and expectations. Further, let
%
%
\begin{equation}\label{eq:tauagain}
\tau^B=\inf\{t\geq0\dvtx B^1_{t} B^2_{t}=0 \}.
\end{equation}
\begin{theorem}\label{thm:theo2}
Let $p > 0$ and $u, v > 0$. Under the above assumptions, the following
conditions are equivalent:
\begin{longlist}
\item
\[
p<\frac{\pi}{{\pi}/{2}+\arctan({\varrho}/({\sqrt
{1-\varrho
^2}}) )},
\]
\item
\[
E^{u,v} [(\tau^B)^{{p}/{2}} ]<\infty,
\]
\item
\[
E^{u,v} [ |(B^1_{\tau^B},B^2_{\tau^B}) |^p ]<\infty.
\]
\end{longlist}
\end{theorem}
\begin{pf}
We start with the proof of the equivalence of (i) and (ii). Define a
cone in the plane with angle $\theta\in(0, 2 \pi)$ by
\[
C(\varphi)= \{re^{i\phi}\dvtx r\geq0, 0\leq\phi\leq\varphi\}
\]
and denote its boundary by $\partial C(\varphi)$. Note that with this
definition, the positive real line is always contained in $C(\varphi)$.
Further, we define, for $\varrho\in(-1,1)$, a sector in $\R^2$ by
\[
S(\varrho) = \biggl\{(x,y)\in\R^2 \dvtx x\geq0, y\geq-\frac{\varrho
}{\sqrt{1-\varrho^2}}x \biggr\}
\]
and denote by $\partial S(\varrho)$ its boundary. Note that this time,
the positive imaginary axis is always in $S(\varrho)$ and that the
angle of the sector at the origin is given by
\[
\theta:=\frac{\pi}{2}+\arctan\biggl(\frac{\varrho}{\sqrt{1-\varrho
^2}} \biggr).
\]
To transform the correlated Brownian motions $B^1,B^2$ to planar
Brownian motion we use the simple fact that $W^1:=B^1,W^2:= (\frac
{B^2-\varrho B^1}{\sqrt{1-\varrho^2}} )$ defines a pair of
independent Brownian motions started in $u, (\frac{v-\varrho
u}{\sqrt{1-\varrho^2}} )$ satisfying $(B^1,B^2)=(W^1,\varrho
W^1+\sqrt{1-\varrho^2}W^2)$. By the definition of $S(\varrho)$, the
planar Brownian motion $(W^1,W^2)$ started in $ (u, (\frac
{v-\varrho u}{\sqrt{1-\varrho^2}} ) )$ hits $\partial S(\varrho)$
if and only if the correlated Brownian motions $B^1,B^2$ started in
$u,v$ hit $\partial C (\frac{\pi}{2} )$. Hence, for $\tau^B$ as
in (\ref{eq:tauagain}), we have
%
%
\begin{equation}\label{tt}
\tau^B=\tau^W:=\inf\{t\geq0 \dvtx(W^1_t,W^2_t) \in\partial
S(\varrho
) \}.
\end{equation}
Since planar Brownian motion is rotation invariant, $S(\varrho)$ may be
rotated to agree with the cone $C(\theta)$,
without changing the exit time. Obviously, with the corresponding
rotated initial conditions, the law of the first exit time $\tau
_{C(\theta)}$ from the cone $C(\theta)$ agrees with the law of $\tau
^W$. For planar Brownian motion in a cone $C(\theta)$ it is well known
(see \cite{S58}, Theorem 2) that
%
%
\begin{equation}
\label{124}
E^{x,y} \bigl[ \bigl(\tau_{C(\theta)} \bigr)^{p/2} \bigr] < \infty
\quad\Leftrightarrow\quad p<\frac{\pi}{\theta},
\end{equation}
independently of $x,y$. (\ref{tt}) and (\ref{124}) now imply the
equivalence of (i) and (ii) and independence of $u,v$.

The proof of the equivalence of (i) and (iii) is via conformal
transformation of the cone $C(\theta)$ to the upper half-plane. Indeed,
we are going to calculate the densities of the exit-point distributions
%
%
\begin{equation}\label{h}
P^{u,v} (B^1_{\tau^B}=0,B^2_{\tau^B}\geq y ),\qquad
P^{u,v} ( B^1_{\tau^B}\leq x,B^2_{\tau^B}=0 ).
\end{equation}
We proceed in three steps: after reducing to independent Brownian
motions in $S(\varrho)$ as for the exit time, we rotate $S(\varrho)$ to
$C(\theta)$ and, finally, stretch the cone to end up with the upper half-plane.

Recall that the first exit of $(B^1,B^2)$ happens at position
$(0,y)\in\partial C(\frac{\pi}{2})$ if and only if the first exit of
$(W^1,W^2)$ takes place at $ (0,\frac{y}{\sqrt{1-\varrho^2}} )\in
\partial S(\varrho)$. Hence, (\ref{h}) transforms to
%
%
\begin{eqnarray}\label{e1}
&&P^{u,v} (B^1_{\tau^B}=0,B^2_{\tau^B}\geq y )
\nonumber\\[-8pt]\\[-8pt]
&&\qquad= P^{u,({v-\varrho u})/{\sqrt{1-\varrho^2}}} \biggl(W^1_{\tau
^W}=0,W^2_{\tau^W}\geq\frac{y}{\sqrt{1-\varrho^2}} \biggr).\nonumber
\end{eqnarray}
In a similar fashion one obtains
%
%
\begin{eqnarray}\label{e2}
&&P^{u,v} ( B^1_{\tau^B}\leq x,B^2_{\tau^B}=0 )\nonumber\\[-8pt]\\[-8pt]
&&\qquad= P^{u,
({v-\varrho u})/{\sqrt{1-\varrho^2}}}
\biggl(W^1_{\tau^W}\leq x,W^2_{\tau^W}=-\frac{\varrho}{\sqrt{1-\varrho
^2}}W^1_{\tau^W} \biggr).\nonumber
\end{eqnarray}
We represent the transformed initial conditions $(z_1,z_2)= (u,\frac
{v-\varrho u}{\sqrt{1-\varrho^2}} )\in S(\varrho)$ in polar
coordinates, that is,
\begin{eqnarray*}
z_1&=&\sqrt{u^2+\frac{(v-\varrho u)^2}{1-\varrho^2}}
\cos\biggl(\arctan\biggl(\frac{v-\varrho u}{u\sqrt{1-\varrho^2}} \biggr)
\biggr),\\ 
z_2&=&\sqrt{u^2+\frac{(v-\varrho u)^2}{1-\varrho^2}}
\sin\biggl(\arctan\biggl(\frac{v-\varrho u}{u\sqrt{1-\varrho^2}} \biggr) \biggr).
\end{eqnarray*}
For the rotation we add the angle $\arctan(\frac{\varrho}{\sqrt
{1-\varrho^2}} )$ to get the new initial condition. Finally, to
map\vspace*{2pt}
the cone $C(\theta)$ conformally to the upper half-plane $\mathbb{H}$,
we apply the map $ z\mapsto z^{\pi/\theta}$ which maps $C(\theta)$
onto $\mathbb{H}$. Using conformal invariance of Brownian motion (e.g.,
Lemma 7.19 of \cite{MP09}), the problem is reduced to the computation
of the exit distribution of planar (time-changed) Brownian motion from
the upper half-plane. Indeed, due to the random time change the (almost
surely finite) exit time changes but not the distribution of the exit
points, which is Cauchy (see Theorem 2.37 of \cite{MP09}). Thus, to
obtain the distribution of the exit points explicitly it, only remains
to specify the transformed initial condition $\tilde{z}_1,\tilde{z}_2$,
which is given by
\begin{eqnarray*}
\tilde{z}_1 &=& \biggl( u^2 +\frac{(v-\varrho u)^2}{1-\varrho^2}
\biggr)^{{\pi}/({2\theta})}\\
&&{}\times\cos\biggl( \frac{\pi}{\theta} \biggl( \arctan
\biggl({\frac{v-\varrho u}{\sqrt{1-\varrho^2}u}} \biggr) +
\arctan\biggl( \frac{\varrho}{\sqrt{1-\varrho^2}} \biggr) \biggr)\biggr),\\%
\tilde{z}_2 &=& \biggl( u^2 +\frac{(v-\varrho u)^2}{1-\varrho^2}
\biggr)^{{\pi}/({2\theta})}\\
&&{}\times\sin\biggl( \frac{\pi}{\theta} \biggl( \arctan
\biggl( {\frac{v-\varrho u}{\sqrt{1-\varrho^2}u}} \biggr) +
\arctan\biggl( \frac{\varrho}{\sqrt{1-\varrho^2}} \biggr) \biggr)\biggr).%
\end{eqnarray*}
Now, let $\tilde{W}^1,\tilde{W}^2$ be two independent Brownian motions
with $\tilde W^1_0=\tilde z_1, \tilde W^2_0=\tilde z_2$ and
\[
\tau^{\tilde W}:= \inf\{ t> 0\dvtx\tilde{W}^2_t=0 \}.
\]
Then, by (\ref{e1}), (\ref{e2}),
\begin{eqnarray*}
P^{u,v} ( B^1_{\tau^B}=0,B^2_{\tau^B}\geq y )&=& P^{u,
({v-\varrho u})/{\sqrt{1-\varrho^2}}} \biggl(W^1_{\tau^W}=0,W^2_{\tau^W}\geq
\frac{y}{\sqrt{1-\varrho^2}} \biggr)\\
&=& P^{\tilde{z}_1,\tilde{z}_2} \biggl(\tilde{W}^1_{\tau^{\tilde W}}\leq
- \biggl(\frac{y}{\sqrt{1-\varrho^2}} \biggr)^{{\pi/\theta}} \biggr),\\
P^{u,v} ( B^1_{\tau^B}\leq x,B^2_{\tau^B}=0 )&=&P^{u,(
{v-\varrho u})/{\sqrt{1-\varrho^2}}} \biggl(W^1_{\tau^W}\leq x,W^2_{\tau
^W}=-\frac{\varrho}{\sqrt{1-\varrho^2}}W^1_{\tau^W} \biggr)\\
&=&P^{\tilde{z}_1,\tilde{z}_2} \biggl( 0\leq\tilde{W}^1_{\tau^{\tilde
W}}\leq\biggl(x \biggl(1+\frac{\varrho^2}{1-\varrho^2} \biggr)^{1/2}
\biggr)^{{\pi/\theta}} \biggr)\\
&=&P^{\tilde{z}_1,\tilde{z}_2} \biggl( 0\leq\tilde{W}^1_{\tau^{\tilde
W}}\leq\biggl(\frac{x}{\sqrt{1-\varrho^2}} \biggr)^{{\pi}/{\theta}} \biggr).
\end{eqnarray*}
Explicit manipulations of the Cauchy distribution yield
%
%
\begin{eqnarray}
\label{h1}
&&P^{u,v} ( B^1_{\tau^B}=0,B^2_{\tau^B}\geq y )\nonumber\\[-8pt]\\[-8pt]
&&\qquad=\int_{y}^{\infty} \frac{1}{\pi\tilde{z}_2\sqrt{1-\varrho
^2}^{\pi
/\theta}}\frac{{\pi}/{\theta}r^{{\pi}/{\theta}-1}
}{1+ (({ ({r}/{\sqrt{1-\varrho^2}} )^{{\pi
}/{\theta}}
+ \tilde{z}_1})/{\tilde{z}_2} )^2 } \,dr,\nonumber\\
\label{h2}
&&P^{u,v} ( B^1_{\tau^B}\leq x ,B^2_{\tau^B}=0 )\nonumber\\[-8pt]\\[-8pt]
&&\qquad=\int_0^x\frac{1}{\pi\tilde{z}_2\sqrt{1-\varrho^2}^{\pi/\theta}}
\frac{ {\pi/\theta} r^{{\pi/\theta}-1}
}{1+ ( ({ ({r}/{\sqrt{1-\varrho^2}} )^{{\pi
/\theta}}-
\tilde{z}_1})/{\tilde{z}_2} )^2} \,dr.\nonumber
\end{eqnarray}
Finally, noting that $\int_0^{\infty}\frac{x^{p+\alpha
-1}}{1+x^{2\alpha
}} \,dx<\infty$ if and only if $p<\alpha$, we deduce from (\ref{h1}) and
(\ref{h2}) that
\[
E^{u,v} [|(B^1_{\tau^B},B^2_{\tau^B})|^p ]<\infty\quad\mbox{if and only if}
\quad p<\frac{\pi}{\theta}.
\]
\upqed\end{pf}

\subsection[Proof of Theorem 2.5]{Proof of Theorem \protect\ref{thm:mc}}\label{subsec:proofthmmc}

The proof relies on a combination of the self-duality based technique
of the proof of Proposition \ref{prop:pb} and the close relation
between the moments of the exit-time and exit-point distribution of
correlated Brownian motions obtained in Theorem \ref{thm:theo2}.
\begin{pf*}{Proof of Theorem \ref{thm:mc}}
We proceed in several steps. First, the result for the nonspatial
model is proved and thereafter the results for the discrete-space and
the continuous-space models. Finally, we present the argument in the
transient case. In the following we use the definition of $B^1,B^2$ and
$\tau$ from Proposition \ref{prop:convlaw}.

\textit{Step} 1. Suppose $(u_t,v_t)$ is a solution of $\mathrm
{SBM}(\varrho,\kappa)_{1,1}$.

``$\Rightarrow$'': We first assume $\varrho<\varrho(p)$, in which case
Theorem \ref{thm:theo2} implies\break $E^{1,1}[\tau^{p/2}]<\infty$. As argued
in the proof of Proposition \ref{prop:sconv}, $u_t$ is a nonnegative
martingale and due to the same arguments satisfies $\E^{1,1}
[[u_{\cdot}]_{t}^{p/2} ]\leq E^{1,1}[\tau^{p/2}]<\infty$ for all
$t\geq0$ and $\kappa>0$. Considering $\bar u_t=u_t-u_0=u_t-1$, we
apply the Burkholder--Davis--Gundy inequality to get
\begin{eqnarray*}
\E^{1,1}[u_t^p]&=&\E^{1,1}[(\bar u_t+1)^p]\\
&=&\E^{1,1}\bigl[\one_{\{\bar u_t\leq1\}}(\bar u_t+1)^p\bigr]+\E^{1,1}\bigl[\one
_{\{
\bar
u_t>1\}}(\bar u_t+1)^p\bigr]\\
&\leq& C_p+C_p\E^{1,1}[{\bar u_t}^p]\\
&\leq& C_p+C_p\E^{1,1}\Bigl[\sup_{0\leq s\leq t} \bar u_s^p\Bigr]\\
&\leq& C_p+C'_p\E^{1,1}[[\bar u_{\cdot}]_t^{p/2}]<\infty
\end{eqnarray*}
independently of $t$ and $\kappa$.

``$\Leftarrow$.'' Conversely, for $\varrho\geq\varrho(p)$, Theorem
\ref{thm:theo2} implies that $E^{1,1}[(B^1_{\tau})^p]=\infty$. Using
Fatou's lemma and almost sure convergence of $u_t$ to $B^1_{\tau}$, the
proof for the nonspatial case is finished with
\[
\liminf_{t\rightarrow\infty}\E^{1,1}[u_t^p]\geq\E^{1,1}[u_{\infty
}^p]=E^{1,1}[(B_{\tau})^p]=\infty.
\]
Again, this lower bound is independent of $\kappa$.

\textit{Step} 2. The proof for $\mathrm{dSBM}({\varrho,\kappa
})_{\one
,\one
}$ is started by reducing the moments for homogeneous initial
conditions to finite initial conditions. Indeed, employing Lem\-ma~\ref
{la:sduality} with $\phi=\psi=\frac{\theta}{2}\mathbf{1}_{k}$, where
$\mathbf{1}_{k}$ denotes the indicator function of site $k \in\Z^d$, gives
\begin{eqnarray*}
\E^{\one,\one} \bigl[ e^{-\sqrt{1-\varrho} \theta(u_t(k)+v_t(k))} \bigr]
&=&\E^{\one,\one} \bigl[e^{-\sqrt{1-\varrho}\langle u_{t}+v_{t},\phi
+\psi
\rangle} \bigr]\\
&=&\E^{\phi,\psi} \bigl[e^{-\sqrt{1-\varrho}\langle\one+\one,\tilde
{u}_{t}+\tilde{v}_{t}\rangle} \bigr]\\
&=&\E^{\one_{k},\one_{k}} \bigl[e^{-\sqrt{1-\varrho}\theta\langle\one,
\tilde
{u}_{t}+\tilde{v}_{t}\rangle} \bigr],
\end{eqnarray*}
where we used the argument of Remark \ref{shift}. Note that, due to our
choice of initial conditions, the complex part of the self-duality
vanishes. Since the above is a Laplace transform identity, we have
\[
\mathcal{L}^{\one,\one} \bigl(u_t(k)+v_t(k) \bigr)=\mathcal{L}^{\mathbf{1}_k,
\mathbf{1}_k} (\langle\one, \tilde{u}_t \rangle+\langle\one,
\tilde{v}_t
\rangle)
\]
and hence
%
%
\begin{equation}\label{23}
\E^{\one,\one} \bigl[ \bigl(u_t(k)+v_t(k)\bigr)^p \bigr] = \E^{\mathbf{1}_k,\mathbf{1}_k}
[ (\langle\one, \tilde{u}_t \rangle+ \langle\one, \tilde{v}_t
\rangle
)^p ].
\end{equation}
We are now prepared to finish the proof of the theorem for the discrete case.

``$\Rightarrow$.'' Suppose $\varrho<\varrho(p)$. Let $M_t=\langle
\one,
\tilde{u}_t \rangle+\langle\one, \tilde{v}_t \rangle$, which due to
Lemma \ref{prop:tmmart} is a square-integrable martingale with
quadratic variation
\[
[M_\cdot]_t = [\langle\one, \tilde{u}_{\cdot}\rangle]_t+
[\langle\one, \tilde{v}_{\cdot}\rangle]_t+2 [\langle\one, \tilde
{u}_{\cdot}
\rangle,\langle\one, \tilde{v}_{\cdot} \rangle]_t=(2+2\varrho)
[\langle\one, \tilde{u}_{\cdot} \rangle]_t.
\]
To apply the Burkholder--Davis--Gundy inequality, we switch again from
$M$ to $\bar M_t=M_t-M_0$, which is a martingale null at zero. Hence,
\[
\E^{\mathbf{1}_k, \mathbf{1}_k} [M_t^p ]=\E^{\mathbf{1}_k, \mathbf{1}_k}
[(\bar{M}_t+M_0)^p ]
\leq C_p+C_p\E^{\mathbf{1}_k, \mathbf{1}_k} [\bar{M}_t^p ].
\]
Then we get from (\ref{23}) and the Burkholder--Davis--Gundy inequality
\begin{eqnarray*}
\E^{\one,\one} \bigl[\bigl(u_t(k)+v_t(k)\bigr)^p \bigr]
&\leq& C_p+C_p\E^{\mathbf{1}_k, \mathbf{1}_k} [\bar{M}_t^p ]\\
&\leq& C_p+C_p\E^{\mathbf{1}_k,\mathbf{1}_k} \Bigl[ \sup_{0\leq s\leq t}
\bar
M_s^p \Bigr]\\
&\leq& C_p+C'_p \E^{\mathbf{1}_k, \mathbf{1}_k} [[\bar M_{\cdot
}]_t^{p/2} ] \\
&=&C_p+ C'_p(2+2\varrho)^{p/2} \E^{\mathbf{1}_k,\mathbf{1}_k}
[[\langle
\one
, \tilde{u}_{\cdot} \rangle]_t^{p/2} ]
\end{eqnarray*}
for some constants $C_p,C'_p$ independent of $t$ and $\kappa$. As in
the proof of Theorem \ref{prop:convlaw}, the random time change which
makes the pair of total masses a pair of correlated Brownian motions is
bounded by $\tau$, that is, $ [\langle\one, \tilde{u}_{\cdot}
\rangle
]_t\leq\tau$ for all $t\geq0$. This yields by Theorem \ref{thm:theo2}
\[
\E^{\one,\one} [u_t(k)^p ]\leq\E^{\one,\one} \bigl[\bigl(u_t(k)+v_t(k)\bigr)^p
\bigr]\leq C_p+C'_p(2+2\varrho)^{p/2} E^{1, 1} [\tau^{p/2} ]<\infty.
\]

``$\Leftarrow$.'' Suppose $\varrho\geq\varrho(p)$. As in the proof of
Theorem \ref{prop:convlaw} we use the almost sure convergence of
$(\langle\one,\tilde{u}_t \rangle,\langle\one, \tilde{v}_t
\rangle)$ to
$(B^1_{\tau}, B^2_{\tau})$. Combining this with Fatou's lemma gives
\begin{eqnarray*}
\liminf_{t\rightarrow\infty} \E^{\mathbf{1}_k,\mathbf{1}_k}
[(\langle
\one, \tilde{u}_t \rangle+ \langle\one, \tilde{v}_t \rangle)^p ]
&\geq&\liminf_{t\rightarrow\infty}\E^{\mathbf{1}_k,\mathbf{1}_k}
[\langle\one, \tilde{u}_t \rangle^p ]\\
&\geq&\E^{\one_k,\one_k} \Bigl[\liminf_{t\rightarrow\infty}\langle
\one,
\tilde{u}_t \rangle^p \Bigr] \\
&=&E^{1,1} [(B^1_{\tau})^p ].
\end{eqnarray*}
The right-hand side is infinite due to Theorem \ref{thm:theo2} and
hence $ \E^{\mathbf{1}_k,\mathbf{1}_k} [(\langle\one, \tilde{u}_t
\rangle+
\langle\one, \tilde{v}_t \rangle)^p ]$ diverges. Equation (\ref{23})
now shows that $\E^{\one,\one}[(u_t(k)+v_t(k))^p]$ also grows
without bound.
Since symbiotic branching processes are nonnegative, this is also
true for
$\E^{\one,\one}[u_t(k)^p]$ as can be seen as follows:
\begin{eqnarray*}
\E^{\one,\one}\bigl[\bigl(u_t(k)+v_t(k)\bigr)^p\bigr]
&\leq&\E^{\one,\one}\bigl[(2u_t(k))^p \one_{\{u_t(k)\geq v_t(k)\}}\bigr]\\
&&{} + \E
^{\one,\one
}\bigl[(2v_t(k))^p \one_{\{u_t(k)<v_t(k)\}}\bigr]\\
&\leq&2^p\E^{\one,\one}[u_t(k)^p]+ 2^p\E^{\one,\one}[v_t(k)^p]\\
&=& 2^{p+1}\E^{\one,\one}[u_t(k)^p],
\end{eqnarray*}
where we used Lemma \ref{la:mdual} to see that $\E^{\one,\one
}[u_t(k)^p]=\E
^{\one,\one}[v_t(k)^p]$.

\textit{Step} 3. The proof for $\mathrm{cSBM}(\varrho,\kappa)_{\one
,\one}$
is slightly more involved since we cannot use the indicator $\one_{x}$ to
get $u_t(x)=\langle u_t,\one_{x}\rangle$, where now $\langle
f,g\rangle
=\int_{\R}f(x)g(x) \,dx$. Instead we use a standard smoothing procedure.
For fixed $x\in\R$ let
\[
p_{\varepsilon}(y)=\frac{1}{\sqrt{2\pi\varepsilon}}e^{-
{(x-y)^2}/({2\varepsilon})},
\]
where we skip the dependence on $x$. The main part is to show that
%
%
\begin{eqnarray}\label{z}
&&\bigl\Vert\bigl(u_t(x)+v_t(x)\bigr)-(\langle u_t,p_{\varepsilon}\rangle+\langle
v_t,p_{\varepsilon}\rangle)\bigr\Vert_{L^p} \nonumber\\[-8pt]\\[-8pt]
&&\qquad\leq\Vert u_t(x)-\langle u_t,p_{\varepsilon}\rangle
\Vert_{L^p}+\Vert v_t(x)-\langle v_t,p_{\varepsilon}\rangle\Vert_{L^p}
\stackrel
{\varepsilon\rightarrow0}{\rightarrow} 0,\nonumber
\end{eqnarray}
which implies
%
%
\begin{equation}\label{z2}
\lim_{\varepsilon\rightarrow0}\Vert\langle u_t,p_{\varepsilon
}\rangle
+\langle
v_t,p_{\varepsilon}\rangle\Vert_{L^p}=\Vert u_t(x)+v_t(x)\Vert_{L^p}.
\end{equation}
Due to symmetry we only consider $\Vert u_t(x)-\langle
u_t,p_{\varepsilon
}\rangle\Vert_{L^p}$. To prove (\ref{z}) we first observe that, due to
the Green function representation provided in Corollary 19 of~\cite{EF04},
\begin{eqnarray*}
&&\Vert u_t(x)-\langle u_t,p_{\varepsilon}\rangle\Vert_{L^p}\\
&&\qquad=
\biggl\Vert P_tu_0(x)-\langle P_{t+\varepsilon},u_0\rangle+\int_0^t\int_{\R
}p_{t-s}(x-b)M(ds,db) \\
&&\qquad\quad\hspace*{85.8pt}{} - \int_0^t\int_{\R}P_{t-s}p_{\varepsilon}(x-b)
M(ds,db) \biggr\Vert_{L^p},
\end{eqnarray*}
where $M(ds,db)$ is a zero-mean martingale measure with quadratic variation
\[
\biggl[\int_0^{\cdot}\int_{\R}f(s,b)M(ds,db) \biggr]_t=\kappa\int_0^t\int
_{\R}f^2(s,b)u_s(b)v_s(b) \,ds\,db
\]
for test functions $f$ such that the integral is well defined (see
Lemma 18 of \cite{EF04} for details).

For homogeneous initial conditions, the first difference vanishes and
it suffices to concentrate on the difference of the stochastic
integrals. By the Burkholder--Davis--Gundy inequality the difference of
the integrals can be estimated as
\begin{eqnarray*}
&&\E\biggl[ \biggl( \int_0^t\int_{\R}p_{t-s}(x-b)M(ds,db)-\int_0^t\int_{\R
}P_{t-s}p_{\varepsilon}(x-b)M(ds,db) \biggr)^{p} \biggr]\\
&&\qquad\leq C\kappa^{p/2}\E\biggl[ \biggl( \int_0^t\int_{\R
}\bigl(p_{t-s}(x-b)-p_{\varepsilon+t-s}(x-b)\bigr)^2u_s(b)v_s(b) \,ds\,db \biggr)^{p/2} \biggr].
\end{eqnarray*}
Now expanding $(p_{t-s}(x-b)-p_{\varepsilon+t-s}(x-b))^2 u_s(b)v_s(b)$ as
\begin{eqnarray*}
&&\bigl(p_{t-s}(x-b)-p_{\varepsilon+t-s}(x-b) \bigr)^{2(p-1)/p}\\
&&\qquad{}\times\bigl(p_{t-s}(x-b)-p_{\varepsilon+t-s}(x-b) \bigr)^{2/p} u_s(b)v_s(b),
\end{eqnarray*}
we get the upper bound (taking the expectation under the integral is
valid since the integrands are nonnegative)
\begin{eqnarray*}
&&C\kappa^{p/2} \biggl[ \biggl(\int_0^t\int_{\R}
\bigl(p_{t-s}(x-b)-p_{\varepsilon+t-s}(x-b) \bigr)^2 \,ds \,db \biggr)^{p-1}\\
&&\qquad\hspace*{10.2pt}{}\times\int_0^t\int_{\R} \bigl(p_{t-s}(x-b)-p_{\varepsilon+t-s}(x-b)
\bigr)^2\E[(u_s(b)v_s(b))^{p} ] \,ds\,db \biggr],
\end{eqnarray*}
where we have used that, for $f, g \in L^p$,
\[
\biggl(\int\bigl(f^{2(p-1)/p}\bigr)(f^{2/p}g) \,dx \biggr)^p\leq\biggl(\int f^2
\,dx \biggr)^{p-1}\int f^2g^p \,dx
\]
by H\"older's inequality. As in \cite{EF04}, page 153, the second term
can now be bounded from above by a constant depending only on $p$ and
$t$. The first factor can be estimated by $\varepsilon^{(p-1)/2}$ due to
\cite{r9}, Lemma 6.2. Hence, for fixed $p>1, x\in\R$ and $t\geq0$,
(\ref{z}) holds and thus we obtain (\ref{z2}).
The rest of the proof is similar to the discrete case but slightly
more technical. Since $p_{\varepsilon}(x-\cdot)$ is rapidly decreasing,
we have
\begin{eqnarray*}
\E^{\one,\one} \bigl[ e^{-2\theta\sqrt{1-\varrho}\langle u_t+v_t,
p_{\varepsilon}\rangle} \bigr]
&=&\E^{\theta p_{\varepsilon},\theta p_{\varepsilon}} \bigl[ e^{-2\sqrt
{1-\varrho}\langle\one,
\tilde{u}_t+\tilde{v}_t \rangle} \bigr]\\
&=&\E^{ p_{\varepsilon}, p_{\varepsilon}} \bigl[e^{-\sqrt{1-\varrho
}2\theta
\langle\one,
\tilde{u}_t+\tilde{v}_t \rangle} \bigr].
\end{eqnarray*}
Thus, we get
\[
\mathcal{L}^{\one,\one} (\langle u_t+v_t,p_{\varepsilon}\rangle
)=\mathcal{L}^{p_{\varepsilon},p_{\varepsilon}} (\langle\one,
\tilde
{u}_t+\tilde{v}_t \rangle)
\]
and in particular
\[
\E^{\one,\one} [(\langle u_t+v_t,p_{\varepsilon}\rangle)^p ]
=\E^{p_{\varepsilon},p_{\varepsilon}} [(\langle\one, \tilde{u}_t
\rangle
+\langle\one, \tilde{v}_t \rangle)^p ].
\]
We may now finish the proof in a similar way to the discrete case.

``$\Rightarrow$.'' Due to (\ref{z2}) we are done if we can bound $\E
^{\one,\one}[\langle u_t+v_t,p_{\varepsilon}\rangle^p]$
independently of
$\varepsilon>0$ and $t\geq0$. This can be done as before: $\langle
\one,
\tilde{u}_t \rangle$ and $\langle\one, \tilde{v}_t \rangle$ are random
time-changed correlated Brownian motions with initial conditions
$\langle\one, p_{\varepsilon} \rangle=1$ for all $\varepsilon>0$.
Using, as
before, the auxiliary martingale
\[
\bar M_t=\langle\one, \tilde{u}_t \rangle+\langle\mathbf
{1},\tilde
{v}_t\rangle-\langle\mathbf{1},\tilde{u}_0 \rangle- \langle
\mathbf{1},
\tilde{v}_0\rangle,
\]
we obtain (as in the discrete case) with the help of the
Burkholder--Davis--Gundy inequality
\begin{eqnarray*}
\E^{\one,\one} \bigl[\bigl(u_t(x)+v_t(x)\bigr)^p \bigr]
&=&\lim_{\varepsilon\rightarrow0}\E^{\one,\one} [\langle
u_t+v_t,p_{\varepsilon}\rangle^{p} ]\\
&=&\lim_{\varepsilon\rightarrow0}
\E^{p_{\varepsilon},p_{\varepsilon}} [\langle\one,\tilde
{u}_t+\tilde{v}_t
\rangle^p ]\\
&\leq& C_p+C_p\lim_{\varepsilon\rightarrow0}
\E^{ p_{\varepsilon},p_{\varepsilon}} [\bar{M}_t^p ]\\
&\leq& C_p+C'_p\lim_{\varepsilon\rightarrow0}
\E^{p_{\varepsilon}, p_{\varepsilon}} [[ \bar{M}_{\cdot
}]^{p/2}_{t} ]\\
&\leq& C_p+ C'_p(2+2\varrho)^{p/2}E^{1,1}[\tau^{p/2}].
\end{eqnarray*}
The positive constants $C_p, C'_p$ are independent of $\varepsilon$ and
$t$, whereas $\bar M$ and the random time change $[\bar M_{\cdot}]_t$
do depend on $\varepsilon$. However, the bound $[\bar{M}_{\cdot
}]_{t}\leq
\tau$ holds for all $\varepsilon>0$ and $t\geq0$ since
$B^1_0=B^2_0=\langle\one,p_{\varepsilon}\rangle=1$. For $\varrho
<\varrho
(p)$ the right-hand side is finite by Theorem \ref{thm:theo2} and
independent of $t\geq0$. Since
$\E^{\one,\one} [u_t(x)^p ]\leq\E^{\one,\one} [(u_t(x)+v_t(x))^p
] $, the first direction is shown.

``$\Leftarrow$.'' First note that by translation invariance of initial
condition, spatial motion and white noise
\[
\E^{\one,\one}\bigl[\bigl(u_t(x)+v_t(x)\bigr)^p\bigr]=\E^{\one,\one}\bigl[\bigl(u_t(y)+v_t(y)\bigr)^p\bigr]
\]
for fixed time $t\geq0$ and arbitrary spatial positions $x,y\in\R$
implying that
\[
\E^{\one,\one}\bigl[\bigl(u_t(x)+v_t(x)\bigr)^p\bigr]=\int_{x-1/2}^{x+1/2}\E^{\one
,\one
}\bigl[\bigl(u_t(y)+v_t(y)\bigr)^p\bigr]\, dy.
\]
Using Fubini's theorem and Jensen's inequality we obtain for $p>1$ the
lower bound
\begin{eqnarray*}
\E^{\one,\one}\bigl[\bigl(u_t(x)+v_t(x)\bigr)^p\bigr]&\geq&\E^{\one,\one} \biggl[ \biggl(\int
_{x-1/2}^{x+1/2}\bigl(u_t(y)+v_t(y)\bigr) \,dy \biggr)^p \biggr]\\
&=&\E^{\one,\one}\bigl[\bigl\langle u_t+v_t,1_{(x-1/2,x+1/2)}\bigr\rangle^p\bigr].
\end{eqnarray*}
We now choose an arbitrary nonnegative (nontrivial) smooth function
$f$ with support contained in $(x-1/2,x+1/2)$ that is bounded by $1$
and integrates to some $c \in(0,1)$, say. A lower bound is now given by
\begin{eqnarray*}
\E^{\one,\one}\bigl[\bigl(u_t(x)+v_t(x)\bigr)^p\bigr]&\geq&\E^{\one,\one}[\langle
u_t+v_t,f\rangle
^p]\\
&=&\E^{f,f}[\langle\tilde u_t+\tilde v_t,\one\rangle^p]\\
&\geq&\E^{f,f}[\langle\tilde u_t,\one\rangle^p],
\end{eqnarray*}
where we utilized for the equality the self-duality of Proposition 5 of
\cite{EF04}.

Finally, as in the discrete case, Fatou's lemma and the martingale
convergence theorem imply
\[
\liminf_{t\to\infty}\E^{\one,\one}\bigl[\bigl(u_t(x)+v_t(x)\bigr)^p\bigr]\geq
E^{c,c}[(B^1_{\tau})^p]=\infty
\]
by Theorem \ref{thm:theo2} and due to nonnegativity of solutions as
well
\[
\liminf_{t\rightarrow\infty}\E^{\one,\one} [u_t(x)^p]=\infty
\]
proving the claim.

\textit{Step} 4. The first direction of the above proof for $\mathrm
{dSBM}(\varrho,\kappa)_{\one,\one}$ also works for the transient
case since
$\E^{\one_k,\one_k} [[\bar{M}_{\cdot}]^{p/2}_{\infty} ]\leq
E^{1,1}[\tau^{p/2}]$ is independent of recurrence/transience.
\end{pf*}

\subsection[Proof of Theorem 2.7]{Proof of Theorem \protect\ref{thm:im}}
We now study the ``criticality'' of the critical curve in more detail.
As a preliminary result (mixed) moments of the nonspatial model are
analyzed. The idea is to combine three different techniques: the
martingale argument which led to Theorem \ref{thm:mc} for $\E
^{1,1}[u_t^n]$, a perturbation argument based on the moment duality
which allows us to deduce exponential increase/decrease of $\E
^{1,1}[u_t^{n-1}v_t]$, and finally moment equations which yield
exponential increase/decrease for all mixed moments $\E^{1,1}[u_t^{n-m}v_t^m]$.
\begin{proposition}\label{0}
The following hold for nonspatial symbiotic branching processes:
\begin{itemize}[(1)]
\item[(1)] For all $\kappa>0$ and $n\in\mathbb{N}$:
\begin{itemize}
\item[$\bullet$] $\E^{1,1}[u_t^n]$ grows to a finite constant if $\varrho
<\varrho(n)$,
\item[$\bullet$] $\E^{1,1}[u_t^n]$ grows subexponentially fast to infinity if
$\varrho=\varrho(n)$,
\item[$\bullet$] $\E^{1,1}[u_t^n]$ grows exponentially fast if $\varrho>\varrho(n)$.
\end{itemize}
\item[(2)] For all $\kappa>0$, $n\in\mathbb{N}$ and $m=1,\ldots,n-1$:
\begin{itemize}
\item[$\bullet$] $\E^{1,1}[u_t^{n-m}v_t^m]$ decreases exponentially fast if
$\varrho<\varrho(n)$,
\item[$\bullet$] $\E^{1,1}[u_t^{n-m}v_t^m]$ neither grows exponentially fast nor
decreases exponentially fast if $\varrho=\varrho(n)$,
\item[$\bullet$] $\E^{1,1}[u_t^{n-m}v_t^m]$ grows exponentially fast if $\varrho
>\varrho(n)$.
\end{itemize}
\end{itemize}
\end{proposition}
\begin{pf}
\textit{Step} 1. Martingale arguments based on the connection of
moments of exit times and exit points of correlated Brownian motions
were carried out in the proof of Theorem \ref{thm:mc}. This led to the
first part of (1). Applying H\"older's inequality with $p=\frac
{n}{n-m}$, $q=\frac{n}{m}$, we get the bound
%
%
\begin{equation}\label{df}
\E^{1,1}[u_t^{n-m}v_t^m]\leq\E^{1,1}[u_t^n]^{(n-m)/n}\E
^{1,1}[v_t^n]^{m/n}=\E^{1,1}[u_t^n]
\end{equation}
by symmetry. This implies that for $\varrho<\varrho(n)$ all mixed
moments stay bounded as well.

\textit{Step} 2. We apply the moment duality for the nonspatial model
as explained in Remark \ref{rdual}. Combining the duality with the
martingale argument of the first step we can understand the case
$\varrho<\varrho(n)$ for mixed moments in a simple way. Note that for
mixed moments the dual process starts with $n-m$ particles of one color
and $m$ particles of the other color at time $0$. Note that for mixed
moments $L_t^{\neq}\geq t$, since there is always at least one pair of
different color. Now suppose $\varrho<\varrho(n)$, then for
$0<\varepsilon
<\varrho(n)-\varrho$ we get
\begin{eqnarray*}
\E^{1,1}[u_t^{n-m}v_t^m]&=&\E\bigl[e^{\kappa(L_t^=+\varrho L_t^{\neq
})}\bigr]=\E
\bigl[e^{\kappa(L_t^=+(\varrho+\varepsilon) L_t^{\neq})}e^{-\kappa
\varepsilon
L_t^{\neq}}\bigr]\\
&\leq&\E\bigl[e^{\kappa(L_t^=+(\varrho+\varepsilon)
L_t^{\neq
})}\bigr]e^{-\kappa\varepsilon t}.
\end{eqnarray*}
Since the first factor of the right-hand side is just the moment $\E
^{1,1}[u_t^{n-m}v_t^m]$ for $\varrho+\varepsilon$ strictly smaller than
$\varrho(n)$, this is bounded for all $t$ and $\kappa$. Hence, for
$\varrho<\varrho(n)$ all mixed moments decrease exponentially fast
proving the first part of (2). Note that since $u_t^n$ is a
submartingale, the moment $\E^{1,1}[u_t^n]$ is nondecreasing.

For $\varrho=\varrho(n)$ we first consider the pure moments. Again,
for the critical case, Theorem \ref{thm:mc} implies
\[
\E\bigl[e^{\kappa(L_t^=+(\varrho(n)-\varepsilon)L_t^{\neq
})}\bigr]<C(\varepsilon
)<\infty
\]
for all $\varepsilon>0$ and $t\geq0$. With the crude estimate
$L_t^{\neq
}\leq{n\choose2}t$ we get
\[
C(\varepsilon)>\E\bigl[e^{\kappa(L_t^=+\varrho(n) L_t^{\neq
})}e^{-\kappa
\varepsilon L_t^{\neq}}\bigr]\geq\E\bigl[e^{\kappa(L_t^=+\varrho(n)L_t^{\neq
})} \bigr]e^{-\kappa\varepsilon{n\choose2}t}.
\]
Since $\varepsilon$ is arbitrary this implies subexponential growth to
infinity of $\E^{1,1} [u_t(k)^n ]$ at the critical point. Hence,
the second part of (1) is proven and combined with (\ref{df}) so is the
upper bound of the second part of (2).

\textit{Step} 3. A direct application of It\^o's lemma and Fubini's
theorem yields
\[
\E^{1,1}[u_t^n]=1+\kappa\pmatrix{n\cr2}\int_0^t\E
^{1,1}[u_s^{n-1}v_s] \,ds.
\]
Since we already know from the martingale arguments that $\E
^{1,1}[u_t^n]$ increases to infinity in the critical case, the mixed
moment $\E[u_t^{n-1}v_t]$ cannot decrease exponentially fast proving
the lower bound of part two of (2). Furthermore, with the same
arguments as above, for $\varrho>\varrho(n)$, this leads to
\[
\E^{1,1}[u_t^{n-1}v_t]=\E\bigl[e^{\kappa(L_t^=+\varrho(n) L_t^{\neq
})}e^{\kappa(\varrho-\varrho(n))L_t^{\neq}} \bigr]\geq\E\bigl[e^{\kappa
(L_t^=+\varrho(n) L_t^{\neq})} \bigr]e^{\kappa(\varrho-\varrho(n))t}.
\]
Since the first factor of the right-hand side equals $\E[u_t^{n-1}v_t]$
at the critical point, it does not decrease exponentially fast. Hence,
the product increases exponentially fast. In particular, due to (\ref
{df}), this also implies the third part of (1). Now it only remains to
prove exponential increase for the other mixed moments. Again, using
It\^o's lemma and Fubini's theorem yields the following moment
equations for the mixed moments:
\begin{eqnarray*}
\E^{1,1}[u_t^{n-2}v_t^2]
&=&1+\kappa\int_0^t\E^{1,1}[u_s^{n-1}v_s] \,ds+\varrho(n-2)\kappa\int
_0^t\E^{1,1}[u_s^{n-2}v_s^2] \,ds\\
&&{} + \pmatrix{{n-2}\cr2 }\kappa\int_0^t\E^{1,1}[u_s^{n-3}v_s^3] \,ds
\end{eqnarray*}
and similarly for all other mixed moments. Since we already know
that\break
$\E^{1,1}[u_t^{n-1}v_t]$ grows exponentially fast in $t$, this implies
exponential growth of $\E^{1,1}[u_t^{n-2}v_t^2]$. Iterating this
argument gives exponential growth of all mixed moments for $\varrho
>\varrho(n)$. This shows the third part of (2) and the proof is finished.
\end{pf}

Now it only remains to prove Theorem \ref{thm:im}, where some ideas for
the nonspatial case are recycled.
\begin{pf*}{Proof of Theorem \ref{thm:im}}
First, due to Lemma \ref{la:mdual}, for homogeneous initial
conditions, the moments of $u_t(k)$ and $v_t(k)$ are equal for all
$t\geq0$. For the existence of the Lyapunov exponents we use a
standard subadditivity argument. Hence, it suffices to show
\[
\E^{\one,\one}[u_{t+s}(k)^n] \leq\E^{\one,\one} [u_{t}(k)^n] \E
^{\one,\one}
[u_{s}(k)^n].
\]
Using Lemma \ref{la:mdual}, we reduce the problem to $\E[e^{\kappa
(L_{t}^=+\varrho L_{t}^{\neq})} ]$, where the dual process $(n_t)$
starts with $n$ particles of the same color all placed at site $k$. By
the tower property and the strong Markov property, we obtain
\[
\E^{n_0} \bigl[e^{\kappa(L_{t+s}^=+\varrho L_{t+s}^{\neq})} \bigr]=\E
^{n_0} \bigl[e^{\kappa(L_{t}^=+\varrho L_{t}^{\neq})}\E^{n_t}
\bigl[e^{\kappa(L_{s}^=+\varrho L_{s}^{\neq})} \bigr] \bigr].
\]
We are done if we can show that
%
%
\begin{equation}
\E^{n'} \bigl[e^{\kappa(L_{s}^=+\varrho L_{s}^{\neq})} \bigr]\leq\E
^{n_0} \bigl[e^{\kappa(L_{s}^=+\varrho L_{s}^{\neq})} \bigr]
\end{equation}
for any given initial configuration $n'$ of the dual process consisting
of $n$ particles. The general initial conditions of the dual process
consist of $n^1$ particles of one color and $n^2$ particles of the
other color ($n^1+n^2=n$) distributed arbitrarily in space at positions
$k_1,\ldots,k_{n}$. Using the duality relation of Lemma \ref{la:mdual},
we obtain
\begin{eqnarray*}
\E^{n'} \bigl[e^{\kappa(L_{s}^=+\varrho L_{s}^{\neq})} \bigr]&=&\E^{\one,\one
}[u_s(k_1)
\cdot\cdot\cdot u_s(k_{n^1})v_s(k_{n^1+1})\cdot\cdot\cdot
v_s(k_{n^1+n^2})]\\
&\leq&\E^{\one,\one} [u_s(k)^{n} ]=\E^{n_0} \bigl[e^{\kappa
(L_{s}^=+\varrho L_{s}^{\neq})} \bigr],
\end{eqnarray*}
where, in the penultimate step, we have used the generalized H\"older
inequality.

Having established existence of the Lyapunov exponents, we now turn to
the more interesting question of positivity. The boundedness for
$\varrho<\varrho(n)$ in Theorem~\ref{thm:mc} immediately implies that
in this case $\gamma(\varrho,\kappa)=0$. Now suppose $\varrho
=\varrho
(n)$, that is, $(\varrho,n)$ lies on critical curve.
We use the perturbation argument which we already used for the
nonspatial case combined with Lemma \ref{la:mdual} and Theorem~\ref
{thm:mc} to prove that in this case moments only grow subexponentially
fast. This implies that the Lyapunov exponents are zero. Again we
switch from $\E^{\one,\one}[u_t(k)^n]$ to $\E[e^{\kappa
(L_t^=+\varrho
L_t^{\neq})} ]$, where the dual process is started with all
particles at the same site and the same color. Since moments below the
critical curve are bounded, we can proceed as for the nonspatial model.
For any $\varepsilon>0$, we get
\begin{eqnarray*}
\infty&>&C(\varepsilon)>\E\bigl[e^{\kappa(L_t^=+\varrho L_t^{\neq
})}e^{-\kappa\varepsilon L_t^{\neq}} \bigr]
\geq\E\bigl[e^{\kappa(L_t^=+\varrho L_t^{\neq})} \bigr]e^{-\kappa
\varepsilon{n\choose2}t}
\\
&\geq&\E^{\one,\one} [u_t(k)^n ]e^{-\kappa\varepsilon{n\choose
2}t},
\end{eqnarray*}
where we estimated the collision time of particles of different colors
by the collision time of all particles which is bounded from above by
${n\choose2}t$. Since $\varepsilon$ on the right-hand side is
arbitrary, $\gamma(\varrho,\kappa)$ cannot be positive.

Finally, we assume $\varrho>\varrho(n)$. The idea is to reduce the
problem to the nonspatial case which we already discussed in
Proposition \ref{0}. Actually, we prove more than stated in the theorem
since we also show that mixed moments $\E^{\one,\one
}[u_t(k)^{n-m}v_t(k)^m]$ grow exponentially fast. For $m=1,\ldots,n-1$ the
perturbation argument leads to
\[
\E^{\one,\one}[u_t(k)^{n-m}v_t(k)^m]=\E\bigl[e^{\kappa(L_t^=+\varrho
L_t^{\neq})} \bigr]=\E\bigl[e^{\kappa(L_t^=+\varrho(n) L_t^{\neq
})}e^{\kappa(\varrho-\varrho(n))L_t^{\neq}} \bigr].
\]
The idea is to obtain a lower bound by conditioning on the event that
all particles have not changed their spatial positions before time $t$
(but, of course, have changed their colors). Under this condition the
particle dual is precisely the particle dual of the nonspatial model.
More precisely, we get the lower bound
\begin{eqnarray*}
&&\E \bigl[e^{\kappa(L_t^=+\varrho(n) L_t^{\neq})}e^{\kappa(\varrho
-\varrho(n))L_t^{\neq}};\mbox{no spatial change of particles before
time } t \bigr]\\
&&\qquad=\E\bigl[e^{\kappa(L_t^=+\varrho(n) L_t^{\neq})}e^{\kappa(\varrho
-\varrho(n))L_t^{\neq}} |\mbox{no spatial change of particles before
time } t \bigr]\\
&&\qquad\quad{} \times\P[\mbox{no spatial change of particles before time
}t ]\\
&&\qquad=\E\bigl[e^{\kappa(L_t^=+\varrho(n) L_t^{\neq})}e^{\kappa(\varrho
-\varrho(n))L_t^{\neq}} |\mbox{no spatial change of particles} \bigr]e^{-nt},
\end{eqnarray*}
where the final equality is valid since the event $\{$no spatial
change of particles before time $t\}$ has probability $e^{-nt}$. This
is true since the event is precisely the event that $n$ independent
exponential clocks with parameter $1$ did not ring before time $t$. For
$1\leq m\leq n-1$ there is always at least one pair of particles of
different colors and, hence, we get the lower bound
\[
\E\bigl[e^{\kappa(L_t^=+\varrho(n) L_t^{\neq})}|\mbox{no spatial change
of particles until time } t
\bigr]e^{\kappa(\varrho-\varrho(n))t}e^{-nt},
\]
which equals
\[
\E^{1,1}[u_t^{n-m}v_t^m]e^{\kappa(\varrho-\varrho(n))t}e^{-nt}
\]
for a nonspatial symbiotic branching process with critical correlation
$\varrho=\varrho(n)$. Choosing $\kappa$ such that $\kappa(\varrho
-\varrho(n))>n$ the result now follows from Proposition~\ref{0}.
\end{pf*}

As mentioned in the course of the proof, we actually proved that for
$\varrho>\varrho(n)$ and $m=0,\ldots,n$
\[
\E^{\one,\one}[u_t(k)^{n-m}(k)v_t^m(k)]
\]
grows exponentially in $t$. As for the nonspatial model one could ask
whether, and if so how fast, mixed moments decrease for $\varrho
<\varrho
(n)$. For the second moments it was shown in \cite{AD09} that for
$\varrho<\varrho(2)=0$
\[
\E^{\one,\one}[u_t(k)v_t(k)]\approx
\cases{
\dfrac{1}{\sqrt{t}}, &\quad $d=1$,\cr
\dfrac{1}{\log(t)}, &\quad $d=2$,\vspace*{2pt}\cr
1, &\quad $d\geq3$,}
\]
where $\approx$ denotes weak asymptotic equivalence as $t\to\infty$. It
would be interesting to see whether or not different rates of decrease
appear for moments.

A detailed quantitative study of the Lyapunov exponents as functions of
$\varrho$ and $\kappa$ has so far only been carried out for second
moments (see \cite{AD09}). In contrast to the parabolic Anderson model,
where higher Lyapunov exponents are well studied (see \cite{GdH07}), we
do not have much insight. Only a first upper bound for the Lyapunov
exponents in $\kappa$ and the distance to the critical curve can be
obtained from the perturbation argument of the previous proof.
\begin{proposition}\label{up}
If $\varrho>\varrho(n)$, then $\gamma_n(\varrho,\kappa)\leq\kappa
{n\choose2}(\varrho-\varrho(n))$.
\end{proposition}
\begin{pf}
By Lemma \ref{la:mdual} and Theorem \ref{thm:mc} for $\varrho
>\varrho
(n)$, there are constants $C(\varepsilon)$ such that
\begin{eqnarray*}
C(\varepsilon)&>&\E\bigl[{e^{\kappa(L_t^=+(\varrho-(\varrho-\varrho
(n))-\varepsilon)L_t^{\neq})}} \bigr]\\
&=&\E\bigl[e^{\kappa(L_t^=+\varrho L_t^{\neq})}e^{-\kappa
(\varrho-\varrho(n)+\varepsilon)L_t^{\neq}} \bigr]\\
&\geq&\E\bigl[e^{\kappa(L_t^=+\varrho L_t^{\neq})} \bigr]e^{-\kappa
(\varrho-\varrho(n)+\varepsilon) {n\choose2}t}.
\end{eqnarray*}
Hence, for all $\varepsilon>0$
\[
\E^{\one,\one}[u_t(k)^n]\leq C(\varepsilon)e^{\kappa(\varrho
-\varrho
(n)+\varepsilon) {n\choose2}t},
\]
yielding the result.
\end{pf}

\section{Speed of propagation of the interface}
\label{sec:wavespeed}
In this section we show how to use the moment bounds of Theorem \ref
{thm:mc} to obtain an improved upper bound on the speed of propagation
of the interface as defined in Definition \ref{def:ifc}. We will only
sketch the crucial parts in the proof of Theorem 6 of \cite{EF04} that
need modification. Note that the method used here is based on Mueller's
``dyadic grid technique'' introduced in \cite{M91}.
\begin{pf*}{Proof of Theorem \ref{cor:wavespeed}}
To prove that the interface will eventually be contained in
\[
\bigl[-C\sqrt{T\log(T)},C\sqrt{T\log(T)} \bigr]
\]
(for suitable $C>0$), by symmetry, it suffices to show that the right
endpoint of the interface
\[
R(u_t):=\sup\{x\in\R| u_t(x)>0 \}
\]
up to time $T$ can eventually be bounded by $C\sqrt{T\log(T)}$. To this
end we define
\[
A_n:= \Bigl\{\sup_{t\leq n}R(u_t)>C\sqrt{n\log(n)} \Bigr\}
\]
and show that, for suitably chosen $C$, $\P^{1_{\R^-},1_{\R^+}}
(\limsup_{n\in\mathbb{N}} A_n )=0$. By the Borel--Cantelli lemma,
this follows from
%
%
\begin{equation}\label{abc}
\sum_{n=0}^{\infty}\P^{1_{\R^-},1_{\R^+}}(A_n)<\infty.
\end{equation}
In the following we modify the arguments of \cite{EF04} to obtain an
upper bound for $\P^{1_{\R^-},1_{\R^+}}(A_n)$ which is sumamble over $n$.
\begin{lemma}\label{la:m}
For any integer $n$ there is a finite constant $c_n$ such that for
$\varrho<\varrho(4n-1)$
\[
\E^{1_{\R^-},1_{\R^+}}[(u_t(x)v_t(x))^n]\leq c_n\sqrt{P_t1_{\R
^-}(x)},\qquad x\in\R, t\geq0.
\]
\end{lemma}
\begin{pf}
First recall from (87) of \cite{EF04} that $\E^{1_{\R^-},1_{\R
^+}}[u_t(x)]=P_t1_{\R^-}(x)$. We now use H\"older's inequality and
Theorem \ref{thm:mc} to reduce the mixed moment to the first moment:
\begin{eqnarray*}
&& \E^{1_{\R^-},1_{\R^+}} [(u_t(x)v_t(x))^n ]\\
&&\qquad=\E^{1_{\R^-},1_{\R^+}} [u_t(x)^{1/2}u_t(x)^{n-1/2}v_t(x)^n ]\\
&&\qquad\leq(\E^{1_{\R^-},1_{\R^+}} [u_t(x) ] )^{1/2} (\E^{\one
,\one} [u_t(x)^{4n-1} ] )^{({2n-1})/({8n-2})}\\
&&\qquad\quad{}\times (\E^{\one,\one}
[v_t(x)^{4n-1} ] )^{{n}/({4n-1})}.
\end{eqnarray*}
This follows from the generalized H\"older inequality with exponents
$2,(8n-2)/(2n-1)$ and $(4n-1)/n$. The first factor yields the heat flow
and Theorem~\ref{thm:mc} shows that the latter two factors are bounded
by constants for $\varrho<\varrho(4n-1)$.
\end{pf}

We now strengthen the estimate of Lemma 23 of \cite{EF04} of the
stochastic part
\[
N_t(b)=\int_0^t\int_{\R}p_{t-s}(b-a)M(ds,da)
\]
of the convolution representation of solutions of Corollary 20 of \cite{EF04}.
\begin{lemma}\label{la:ma}
For $\varrho<\varrho(35)$ there is a constant $C_3$ such that for
$\varepsilon\in(0,1)$, $A,T\geq1$, the following estimate holds:
\[
\P^{1_{\R^-},1_{\R^+}} \bigl(|N_t(b)|\geq\varepsilon\mbox{ for some }
t\leq T\mbox{ and }b\geq A \bigr)\leq C_3 \varepsilon^{-18}\frac
{T^{22}}{\sqrt{A}}p_{2T}(A).
\]
\end{lemma}
\begin{pf}
The proof is along the same lines of \cite{EF04} replacing only in
(116) the weaker (exponentially growing) moment bound of \cite{EF04} by
our stronger (bounded) moment bound. In the following we sketch the
arguments to show where the moments appear. Before performing the
``dyadic grid technique,'' increments of $N_t$ need to be estimated.
First, by definition
\begin{eqnarray*}
&&\E^{1_{\R^-},1_{\R^+}}[|N_t(a)-N_{t'}(a')|^{2q}]\\
&&\qquad=\E^{1_{\R
^-},1_{\R
^+}} \biggl[ \biggl|\int_0^t\int_{\R}\bigl(p_{t-s}(b-a)-p_{t'-s}(b-a')\bigr)M(ds,db)
\biggr|^{2q} \biggr],
\end{eqnarray*}
which by Burkholder--Davis--Gundy and H\"older's inequality gives the
upper bound
\begin{eqnarray*}
&&C_1 \biggl|\int_0^t\int_{\R}[p_{t-s}(b-a)-p_{t'-s}(b-a')]^2 \,db \,ds
\biggr|^{q-1}\\
&&\qquad{}\times\int_0^t\int_{\R}[p_{t-s}(b-a)-p_{t'-s}(b-a')]^2 \E^{1_{\R
^-},1_{\R^+}} [(u_s(b)v_s(b))^q ] \,db \,ds.
\end{eqnarray*}

Using Lemma \ref{la:m} and classical heat kernel estimates we can
derive (see the calculation on pages 153, 154 of \cite{EF04}) the upper bound
\begin{eqnarray*}
&& \E^{1_{\R^-},1_{\R^+}}[|N_t(a)-N_{t'}(a')|^{2q}]\\
&&\qquad\leq C_2\bigl((|t'-t|^{1/2}+|a'-a|)\wedge t^{1/2}\bigr)^{q-1} \bigl(\sqrt
{tP_t1_{\R^-}(a)}+\sqrt{t'P_{t'}1_{\R^-}(a')} \bigr).
\end{eqnarray*}
This upper bound corresponds to (119) of \cite{EF04} where they have an
additional exponentially growing factor coming from their moment bound.
The dyadic grid technique can now be carried out as in \cite{EF04},
choosing $q=9$, without carrying along their exponential factor. Hence,
we may delete the exponential term from their final estimate (110).
Note that the necessity of $\varrho<\varrho(35)$ comes from our choice
$q=9$ and Lemma \ref{la:m}.
\end{pf}

The following lemma corresponds to Proposition 24 of \cite{EF04}.
\begin{lemma}
If $\varrho<\varrho(35)$ then, for some constants $C_4, C_5$, the
following estimate holds for $T\geq1$ and $r\geq C_4 \sqrt{T}$:
\[
\P^{1_{\R^-},1_{\R^+}} \Bigl(\sup_{t\leq T}R(u_t)>r \Bigr)\leq C_5 T^{22}
p_{16T}(r).
\]
\end{lemma}
\begin{pf}
All we need to do is to argue that Proposition 24 of \cite{EF04} is
valid for $r\geq C_4 \sqrt{T}$ instead of $r\geq9^4(1\vee\kappa) T$.
We perform the same decomposition and note that the estimates of Step 2
of \cite{EF04} are already given for $r\geq C_4 \sqrt{T}$ if $C_4$ is
large enough. The only trouble occurs in their Step 3. Up to the
estimate (154), this step works for $r\geq C_4 \sqrt{T}$ but here their
(weaker) Lemma 23 produces an exponential in $T$. More precisely, they
need to justify
\[
e^{9^5\kappa^2 T/c}\frac{T^{22}}{\sqrt{r}}p_{8T}(r)\leq
T^{22}p_{16T}(r),
\]
which is only valid for $r\geq9^4 (1\vee\kappa)T$. As our Lemma \ref
{la:ma} avoids the exponential on the left-hand side the estimate holds
for $r\geq C_4\sqrt{T}$ with suitably chosen $C_4$ and $C_5$.
\end{pf}

The significant distinction of the previous lemma to the result of
\cite{EF04} is that the inequality is not only valid for $r\geq
9^4(1\vee\kappa)T$ but for $r\geq C_4\sqrt{T}$. At this point one
might hope to obtain a square-root upper bound for the growth of the
interface but this fails in the final step in which we validate (\ref{abc}):
\begin{eqnarray*}
\sum_{n=0}^{\infty}\P^{1_{\R^-},1_{\R^+}}(A_n)&\leq&\sum
_{n=0}^{\infty
}C_5 n^{22} p_{16n} \bigl(C\sqrt{n\log(n)} \bigr)\\
&=&\sum_{n=0}^{\infty}C_5 n^{22} \frac{1}{\sqrt{\pi32n}}e^{-
{C^2n\log(n)}/({32n})}\\
&=&\frac{C_5}{\sqrt{32\pi}}\sum_{n=0}^{\infty}n^{22-C^2/32-1/2},
\end{eqnarray*}
which is finite for $C$ large enough.
\end{pf*}

\section*{Acknowledgments}
This work is part of the Ph.D. thesis of the second author who would
like to thank the students and faculty from TU Berlin for many
discussions.

The authors would like to express their gratitude to an anonymous
referee for a very careful reading of the manuscript and for pointing
out and correcting an error in an earlier version of Theorem
\ref{cor:wavespeed}.

%

%
\printaddresses


\begin{thebibliography}{25}

\bibitem{AD09}
%
\begin{barticle}[vtex]
\bauthor{\bsnm{Aurzada},~\bfnm{F.}\binits{F.}} \AND
\bauthor{\bsnm{D\"oring},~\bfnm{L.}\binits{L.}}
(\byear{2010}).
\btitle{Intermittency and aging for the symbiotic branching
model}.
\bjournal{Ann. Inst. H. Poincar\'e}.
\bnote{To appear}.
\end{barticle}
%
\endbibitem

\bibitem{CDG04}
%
\begin{barticle}[mr]
\bauthor{\bsnm{Cox},~\bfnm{J.~T.}\binits{J.~T.}},
\bauthor{\bsnm{Dawson},~\bfnm{D.~A.}\binits{D.~A.}} \AND
\bauthor{\bsnm{Greven},~\bfnm{A.}\binits{A.}}
(\byear{2004}).
\btitle{Mutually catalytic super branching random walks: Large finite systems
and renormalization analysis}.
\bjournal{Mem. Amer. Math. Soc.}
\bvolume{171}
\bpages{viii+97}.
\bid{mr={2074427}}
\end{barticle}
%
\endbibitem

\bibitem{CK00}
%
\begin{barticle}[mr]
\bauthor{\bsnm{Cox},~\bfnm{J.~Theodore}\binits{J.~T.}} \AND
\bauthor{\bsnm{Klenke},~\bfnm{Achim}\binits{A.}}
(\byear{2000}).
\btitle{Recurrence and ergodicity of interacting particle systems}.
\bjournal{Probab. Theory Related Fields}
\bvolume{116}
\bpages{239--255}.
\bid{doi={10.1007/PL00008728}, mr={1743771}}
\end{barticle}
%
\endbibitem

\bibitem{CKP00}
%
\begin{bincollection}[mr]
\bauthor{\bsnm{Cox},~\bfnm{J.~Theodore}\binits{J.~T.}},
\bauthor{\bsnm{Klenke},~\bfnm{Achim}\binits{A.}} \AND
\bauthor{\bsnm{Perkins},~\bfnm{Edwin~A.}\binits{E.~A.}}
(\byear{2000}).
\btitle{Convergence to equilibrium and linear systems duality}.
In \bbooktitle{Stochastic Models ({O}ttawa, {ON}, 1998)}.
\bseries{CMS Conf. Proc.}
\bvolume{26}
\bpages{41--66}.
\bpublisher{Amer. Math. Soc.}, \baddress{Providence, RI}.
\bid{mr={1765002}}
\end{bincollection}
%
\endbibitem

\bibitem{CM94}
%
\begin{barticle}[mr]
\bauthor{\bsnm{Carmona},~\bfnm{Ren{\'e}~A.}\binits{R.~A.}} \AND
\bauthor{\bsnm{Molchanov},~\bfnm{S.~A.}\binits{S.~A.}}
(\byear{1994}).
\btitle{Parabolic {A}nderson problem and intermittency}.
\bjournal{Mem. Amer. Math. Soc.}
\bvolume{108}
\bpages{viii+125}.
\bid{mr={1185878}}
\end{barticle}
%
\endbibitem

\bibitem{DFX05}
%
\begin{barticle}[mr]
\bauthor{\bsnm{Dawson},~\bfnm{Donald~A.}\binits{D.~A.}},
\bauthor{\bsnm{Fleischmann},~\bfnm{Klaus}\binits{K.}} \AND
\bauthor{\bsnm{Xiong},~\bfnm{Jie}\binits{J.}}
(\byear{2005}).
\btitle{Strong uniqueness for cyclically symbiotic branching diffusions}.
\bjournal{Statist. Probab. Lett.}
\bvolume{73}
\bpages{251--257}.
\bid{doi={10.1016/j.spl.2005.03.012}, mr={2179284}}
\end{barticle}
%
\endbibitem

\bibitem{DP98}
%
\begin{barticle}[mr]
\bauthor{\bsnm{Dawson},~\bfnm{Donald~A.}\binits{D.~A.}} \AND
\bauthor{\bsnm{Perkins},~\bfnm{Edwin~A.}\binits{E.~A.}}
(\byear{1998}).
\btitle{Long-time behavior and coexistence in a mutually catalytic branching
model}.
\bjournal{Ann. Probab.}
\bvolume{26}
\bpages{1088--1138}.
\bid{doi={10.1214/aop/1022855746}, mr={1634416}}
\end{barticle}
%
\endbibitem

\bibitem{EF04}
%
\begin{barticle}[mr]
\bauthor{\bsnm{Etheridge},~\bfnm{Alison~M.}\binits{A.~M.}} \AND
\bauthor{\bsnm{Fleischmann},~\bfnm{Klaus}\binits{K.}}
(\byear{2004}).
\btitle{Compact interface property for symbiotic branching}.
\bjournal{Stochastic Process. Appl.}
\bvolume{114}
\bpages{127--160}.
\bid{doi={10.1016/j.spa.2004.05.006}, mr={2094150}}
\end{barticle}
%
\endbibitem

\bibitem{GdH07}
%
\begin{barticle}[mr]
\bauthor{\bsnm{Greven},~\bfnm{A.}\binits{A.}} \AND\bauthor
{\bparticle{den
}\bsnm{Hollander},~\bfnm{F.}\binits{F.}}
(\byear{2007}).
\btitle{Phase transitions for the long-time behavior of interacting
diffusions}.
\bjournal{Ann. Probab.}
\bvolume{35}
\bpages{1250--1306}.
\bid{doi={10.1214/009117906000001060}, mr={2330971}}
\end{barticle}
%
\endbibitem

\bibitem{GM90}
%
\begin{barticle}[mr]
\bauthor{\bsnm{G{\"a}rtner},~\bfnm{J.}\binits{J.}} \AND
\bauthor{\bsnm{Molchanov},~\bfnm{S.~A.}\binits{S.~A.}}
(\byear{1990}).
\btitle{Parabolic problems for the {A}nderson model. {I}.
{I}ntermittency and
related topics}.
\bjournal{Comm. Math. Phys.}
\bvolume{132}
\bpages{613--655}.
\bid{mr={1069840}}
\end{barticle}
%
\endbibitem

\bibitem{KM09a}
%
\begin{bmisc}[vtex]
\bauthor{\bsnm{Klenke},~\bfnm{A.}\binits{A.}} \AND
\bauthor{\bsnm{Mytnik},~\bfnm{L.}\binits{L.}}
(\byear{2009}).
\bhowpublished{Infinite rate mutually catalytic branching in infinitely
many colonies: Construction, characterization and convergence. Preprint. Available at}
\href{http://arxiv.org/abs/arXiv:0901.0623v1}{arXiv:0901.0623v1}.
\end{bmisc}
%
\endbibitem

\bibitem{KM09b}
%
\begin{bmisc}[vtex]
\bauthor{\bsnm{Klenke},~\bfnm{A.}\binits{A.}} \AND
\bauthor{\bsnm{Mytnik},~\bfnm{L.}\binits{L.}}
(\byear{2009}).
\bhowpublished{Infinite rate mutually catalytic branching in infinitely
many colonies: The longtime behaviour. Preprint. Available at}
\href{http://arxiv.org/abs/arXiv:0901.4120v1}{arXiv:0901.4120v1}.
\end{bmisc}
%
\endbibitem

\bibitem{KO09}
%
\begin{barticle}[mr]
\bauthor{\bsnm{Klenke},~\bfnm{Achim}\binits{A.}} \AND
\bauthor{\bsnm{Oeler},~\bfnm{Mario}\binits{M.}}
(\byear{2010}).
\btitle{A {T}rotter-type approach to infinite rate mutually catalytic
branching}.
\bjournal{Ann. Probab.}
\bvolume{38}
\bpages{479--497}.
\bid{doi={10.1214/09-AOP488}, mr={2642883}}
\end{barticle}
%
\endbibitem

\bibitem{KS98}
%
\begin{bbook}[mr]
\bauthor{\bsnm{Karatzas},~\bfnm{Ioannis}\binits{I.}} \AND
\bauthor{\bsnm{Shreve},~\bfnm{Steven~E.}\binits{S.~E.}}
(\byear{1991}).
\btitle{Brownian Motion and Stochastic Calculus},
\bedition{2nd} ed.
\bseries{Graduate Texts in Mathematics}
\bvolume{113}.
\bpublisher{Springer}, \baddress{New York}.
\bid{mr={1121940}}
\end{bbook}
%
\endbibitem

\bibitem{MP09}
%
\begin{bbook}[vtex]
\bauthor{\bsnm{M{\"o}rters},~\bfnm{Peter}\binits{P.}} \AND
\bauthor{\bsnm{Peres},~\bfnm{Yuval}\binits{Y.}}
(\byear{2010}).
\btitle{Brownian Motion}.
\bpublisher{Cambridge Univ. Press}, \baddress{Cambridge}.
\bid{mr={2604525}}
\end{bbook}
%
\endbibitem

\bibitem{MT97}
%
\begin{barticle}[vtex]
\bauthor{\bsnm{Mueller},~\bfnm{C.}\binits{C.}} \AND
\bauthor{\bsnm{Tribe},~\bfnm{R.}\binits{R.}}
(\byear{1997}).
\btitle{Finite width for a random stationary interface}.
\bjournal{Electron. J. Probab.}
\bvolume{2}
\bpages{1--27}.
\bid{mr={1485116}}
\end{barticle}
%
\endbibitem

\bibitem{M91}
%
\begin{barticle}[mr]
\bauthor{\bsnm{Mueller},~\bfnm{Carl}\binits{C.}}
(\byear{1991}).
\btitle{On the support of solutions to the heat equation with noise}.
\bjournal{Stochastics Stochastics Rep.}
\bvolume{37}
\bpages{225--245}.
\bid{mr={1149348}}
\end{barticle}
%
\endbibitem

\bibitem{M99}
%
\begin{barticle}[mr]
\bauthor{\bsnm{Mytnik},~\bfnm{Leonid}\binits{L.}}
(\byear{1998}).
\btitle{Uniqueness for a mutually catalytic branching model}.
\bjournal{Probab. Theory Related Fields}
\bvolume{112}
\bpages{245--253}.
\bid{doi={10.1007/s004400050189}, mr={1653845}}
\end{barticle}
%
\endbibitem

\bibitem{R95}
%
\begin{bincollection}[mr]
\bauthor{\bsnm{Rebholz},~\bfnm{Joachim~A.}\binits{J.~A.}}
(\byear{1995}).
\btitle{A skew-product representation for the generator of a two sex population
model}.
In \bbooktitle{Stochastic Partial Differential Equations ({E}dinburgh, 1994)}.
\bseries{London Mathematical Society Lecture Note Series}
\bvolume{216}
\bpages{230--240}.
\bpublisher{Cambridge Univ. Press}, \baddress{Cambridge}.
\bid{mr={1352745}}
\end{bincollection}
%
\endbibitem

\bibitem{S80}
%
\begin{barticle}[mr]
\bauthor{\bsnm{Shiga},~\bfnm{Tokuzo}\binits{T.}}
(\byear{1980}).
\btitle{An interacting system in population genetics}.
\bjournal{J. Math. Kyoto Univ.}
\bvolume{20}
\bpages{213--242}.
\bid{mr={582165}}
\end{barticle}
%
\endbibitem

\bibitem{r7}
%
\begin{barticle}[mr]
\bauthor{\bsnm{Shiga},~\bfnm{Tokuzo}\binits{T.}}
(\byear{1992}).
\btitle{Ergodic theorems and exponential decay of sample paths for certain
interacting diffusion systems}.
\bjournal{Osaka J. Math.}
\bvolume{29}
\bpages{789--807}.
\bid{mr={1192741}}
\end{barticle}
%
\endbibitem

\bibitem{r9}
%
\begin{barticle}[mr]
\bauthor{\bsnm{Shiga},~\bfnm{Tokuzo}\binits{T.}}
(\byear{1994}).
\btitle{Two contrasting properties of solutions for one-dimensional stochastic
partial differential equations}.
\bjournal{Canad. J. Math.}
\bvolume{46}
\bpages{415--437}.
\bid{mr={1271224}}
\end{barticle}
%
\endbibitem

\bibitem{S58}
%
\begin{barticle}[mr]
\bauthor{\bsnm{Spitzer},~\bfnm{Frank}\binits{F.}}
(\byear{1958}).
\btitle{Some theorems concerning {$2$}-dimensional {B}rownian motion}.
\bjournal{Trans. Amer. Math. Soc.}
\bvolume{87}
\bpages{187--197}.
\bid{mr={0104296}}
\end{barticle}
%
\endbibitem

\bibitem{T95}
%
\begin{barticle}[vtex]
\bauthor{\bsnm{Tribe},~\bfnm{Roger}\binits{R.}}
(\byear{1995}).
\btitle{Large time behavior of interface solutions to the heat
equation with
{F}isher--{W}right white noise}.
\bjournal{Probab. Theory Related Fields}
\bvolume{102}
\bpages{289--311}.
\bid{doi={10.1007/BF01192463}, mr={1339735}}
\end{barticle}
%
\endbibitem

\bibitem{W86}
%
\begin{bincollection}[mr]
\bauthor{\bsnm{Walsh},~\bfnm{John~B.}\binits{J.~B.}}
(\byear{1986}).
\btitle{An introduction to stochastic partial differential equations}.
In \bbooktitle{\'{E}cole D'\'et\'e de Probabilit\'es de {S}aint-{F}lour,
{XIV}---1984}.
\bseries{Lecture Notes in Math.}
\bvolume{1180}
\bpages{265--439}.
\bpublisher{Springer}, \baddress{Berlin}.
\bid{mr={876085}}
\end{bincollection}
%
\endbibitem

\end{thebibliography}
\end{document}